\documentclass[xypic,amscd,syntonly,amssymb]{amsart}

\usepackage{graphicx}
\usepackage[all]{xy}
\usepackage{amssymb}
\usepackage{amsmath}
\usepackage{epsfig}
\theoremstyle{plain}
\newtheorem{Thm}{Theorem}
\newtheorem{Prop}[Thm]{Proposition}
\newtheorem{Cor}[Thm]{Corollary}
\newtheorem{Lem}[Thm]{Lemma}

 \theoremstyle{definition}
\theoremstyle{remark}

\numberwithin{equation}{section}

\begin{document}
 \title{Lifting Hamiltonian loops to isotopies in fibrations}

 \author{ ANDR\'{E}S   VI\~{N}A}
\address{Departamento de F\'{i}sica. Universidad de Oviedo.   Avda Calvo
 Sotelo.     33007 Oviedo. Spain. }
 \email{vina@uniovi.es}
\thanks{This work has been partially supported by Ministerio de Educaci\'on y
 Ciencia, grant   FPA2009-11061}

\maketitle

\begin{abstract}
Let $G$ be a Lie group, $H$ a closed subgroup and $M$ the
homogeneous space $G/H$. Each representation $\Psi$ of $H$
determines a $G$-equivariant principal bundle ${\mathcal P}$ on
$M$ endowed with a $G$-invariant connection. We consider subgroups
${\mathcal G}$ of the diffeomorphism group ${\rm Diff}(M)$, such
that, each vector field $Z\in{\rm Lie}({\mathcal G})$ admits a
lift to a preserving connection  vector field on ${\mathcal P}$.
We  prove that $\#\,\pi_1({\mathcal
G})\geq \#\,\Psi(Z(G))$. This relation is applicable to
  subgroups ${\mathcal G}$ of the Hamiltonian groups of
the flag varieties of a semisimple group $G$.
%, equipped with a certain symplectic structures.

Let $M_{\Delta}$ be the toric manifold determined by the Delzant
polytope $\Delta$. We put $\varphi_{\bf b}$ for the the loop in
the Hamiltonian group of $M_{\Delta}$ defined by the lattice
vector ${\bf b}$. We give a sufficient condition, in terms of the
mass center of $\Delta$, for the loops $\varphi_{\bf b}$ and
$\varphi_{\bf\tilde b}$ to be homotopically inequivalent.

\end{abstract}
   \smallskip

%\subjclass[2000]{ 53D50, 22E45}


%%%%%%%%%%%%%%%%%%%%%%%%%%%%%%%%%%%%%%%%%%%%%%%%%%%%%%%%%%%%%%%%%%%%%%%%%%%%%%%%%%%%%%%%%%%%%%%%%%%%%%%%%%%%%%
%%%%%%%%%%%%%%%%%%%%%%%%%%%%%%%%%%%%%%%%%%%%%%%%%%%%%%%%%%%%%%%%%%%%%%%%%%%%%%%%%%%%%%%%%%%%%%%%%%%%%%%%%%%%%%%%%%%%

\section{Introduction}
In this note we will concern with the homotopy of some Lie
subgroups ${\mathcal G}$ of ${\rm Diff}(M)$, the diffeomorphism
group of a manifold $M$ in two cases: (I) When $M$ is  a
homogeneous space. (II) When $M$ is a symplectic toric manifold.

As it is well-known, the lifting of loops in a space $X$ to paths
in a fibration over $X$ is sometimes used for
  studying the
  %the study of
homotopy of $X$. Here, we  will also lift loops in ${\mathcal G}$
to isotopies in appropriate principal fiber  bundles   on $M$, in
order to deduce consequences about $\pi_1({\mathcal G})$.
% We will follow this method in the both mentioned cases.

\smallskip
\noindent
{\bf Case (I).}
  Given a connected Lie group $G$.  Let $(H,\Psi)$ be a pair consisting of
   a closed subgroup $H$ of $G$  and
  a representation $\Psi$  of $H$. By means the representation
   and a
 complement ${\mathfrak l}$ of ${\mathfrak h}$ in ${\mathfrak
g}$ (see (\ref{frakl0})),
  we will
construct a $G$-equivariant principal fibre bundle ${\mathcal P}$
over $M=G/H$ and a $G$-invariant connection $\Omega$ on it.

 We will consider groups ${\mathcal G}$ of ${\rm Diff}(M)$ such that each vector field
 %The Lie algebras of the groups ${\mathcal G}$ we will consider
 %satisfy the following property:
 %Each vector field
  $Z\in{\rm Lie}({\mathcal G})$ admits a lift $U$ to a vector field on
${\mathcal P}$ with the following properties: $U$ is invariant
under the natural right action on ${\mathcal P}$ and it is an
infinitesimal symmetry of $\Omega$ (the usual horizontal lift does
not satisfy the last property). These properties allow us to lift
a loop $\psi=\{\psi_t\}_{t\in[0,\,1]}$ in the group ${\mathcal G}$
to an isotopy  $\{F_t\}_{t\in[0,\,1]}$ in ${\mathcal P}$, whose
final element
%. Moreover, $F_1$ turns out to be with
$F_1$  is a preserving connection gauge transformation of
${\mathcal P}$. The comparison of the gauge transformations,
obtained by lifting different loops in ${\mathcal G}$, will permit
us to distinguish some non-homotopic loops. For the sake of
brevity, we will say that a subalgebra ${\mathfrak X}$ of
${\mathfrak X}(M)$ admits a lift to ${\mathcal P}$ if its elements
admit  lifts to vector fields on ${\mathcal P}$ with the mentioned
properties.

 If $\Psi$ is a
representation of $H$ on the vector space $W$, then the above
mentioned bundle ${\mathcal P}$ is a ${\rm GL}(W)$-principal
bundle and the gauge transformation $F_1$ associated to a loop
$\psi$ in ${\mathcal G}$ is defined by a map
 \begin{equation}\label{maph}
h(\psi):{\mathcal
 P}\to{\rm GL}(W).
 \end{equation}
We will prove    that $h(\psi)$ is constant on each   orbit
 $\{F_t(p)\,|\,t\in[0,\,1]\}$ ({\bf Proposition \ref{fconstant}}). This fact will be important for
characterizing the homotopy class of some loops in ${\mathcal G}$.

In the case that the left multiplication of points in $M=G/H$ by
each $g\in G$   is an element of ${\mathcal G}$ and $Z(G)$ -the
center of $G$- is a subgroup of $H$, then a path
$\{g_t\}_{t\in[0,\,1]}$
 in $G$ with initial point at the identity element $e$ of $G$ and final point $g_1$ in $Z(G)$
 determines a loop $\varphi$ in ${\mathcal G}$.
These loops $\varphi$ enjoy nice properties:

(a) The corresponding gauge function $h(\varphi)$ is constant and
equal to
 %simply the multiplication by the constant
$\Psi(g_1)\in{\rm GL}(W)$ ({\bf Proposition \ref{PropH10}}).
%We will denote by ${\mathcal L}$ the set of loops in ${\mathcal
%G}$ defined for such paths.

 (b) Under certain hypotheses, we will
 construct for each $g_1\in Z(G)$ a particular path $\{g_t\}$ in  $G$, such that
for the corresponding loop $\varphi$ in ${\mathcal G}$ the
following property  holds: If $\{\varphi^s\}_s$ is a deformation
of the loop $\varphi$ in  ${\mathcal G}$, then for all $s$ the map
$h(\varphi^s)$ is the constant map that takes the value
$\Psi(g_1)$. Thus, $\#\,\pi_1({\mathcal G})\geq\#\,\Psi(Z(G))$
(see {\bf Theorem \ref{pi1G}}).

%Let ${\bf G}$ denote the group consisting of the diffeomorphisms
%of $M$ defined by the left multiplication by the elements of $G$.
%As a corollary  of Theorem \ref{pi1G},  we will deduce that
%$\#\,\pi_1({\bf G})\geq \#\,\Psi(Z(G))$, for any irreducible
%representation $\Psi$ of $H$ ({\bf Corollary \ref{CorbfG}}).

%As we will show, the number $\#\,\Psi(Z(G))$ is a lower bound for
%the cardinal of the fundamental group of several subgroups of
%${\rm Diff}(M)$ related with the representation $\Psi$.

\smallskip

 The case when $\Psi$ is
a {\em $1$-dimensional} representation  presents three related
properties:

(i) The curvature ${\bf K}$ of the connection ${\Omega}$ on
${\mathcal P}$ projects a closed $2$-form $\omega$ on $M$; thus,
$\omega$ defines a presymplectic structure on $M$. (Obviously, the
form $\omega$ can also be defined directly from
 the $1$-dimensional representation $\Psi$ of
$H$ and the complement ${\mathfrak l}$ of ${\mathfrak h}$).

  (ii) For
any loop $\psi$ in a group ${\mathcal G}$ whose Lie algebra admits
a lift to ${\mathcal P}$, the map (\ref{maph}) is constant, i.e.,
the gauge transformation
 $F_1$ determined by $\psi$ is simply the multiplication by a constant
$\Theta(\psi)$
 %that will be denoted $\Theta(\psi)$
 (see {\bf Proposition \ref{fconstant2}}).

 (iii) For each vector field $Z$ belonging to ${\rm Lie}({\mathcal G})$,
   there exists a function $J_Z$ on $M$ such that  ${\rm d}J_Z=-\iota_Z\omega.$
   That is,  $J_Z$ is  a Hamiltonian
function for $Z$ relative to $\omega$.

If the first homology group $H_1(M,\, {\mathbb Z})$ vanishes, then
$\Theta(\psi)$ is an ``action
 integral" (relative to the presymplectic structure $\omega$) around the loop $\psi$
(see (\ref{Theta})).
 However,
since the functions $J_Z$ will be not ``normalized",  this
``action
 integral" is  not invariant under deformations of the loop $\psi$ in ${\mathcal G}$.

We will consider the algebra ${\mathfrak X}_{\omega}$ consisting
of all the vector fields $Z$ on $M$, such that $\iota_Z\omega$ is
an exact $1$-form; i.e.,
%{\mathfrak
%L}_Z\omega=0$ (${\mathfrak L}=$ Lie derivative); i. e.,
 the $\omega$-Hamiltonian vector fields on $M$.  The   respective connected
Lie group is denoted by ${\mathcal H}_{\omega}(M).$
 If it is
possible to choose for each $Z\in{\mathfrak X}_{\omega}$ an
$\omega$-Hamiltonian function $J_Z$, so that,  $J_Z$ depends
continuously on   $Z$,  then the algebra ${\mathfrak X}_{\omega}$
admits a lift to ${\mathcal P}$.
  In this
case,  we will prove that
$\#\,\pi_1({\mathcal H}_{\omega}(M))\geq\#\Psi(Z(G))$ ({\bf
Theorem \ref{PropHomega}}).
 %From other point of view, this theorem is a consequence of the
% fact that now
%$\Theta(\varphi)$ is invariant under
%deformations of $\varphi$ in ${\mathcal H}_{\omega}(M).$

\smallskip

When $G$ is a linear complex semisimple Lie group,
 %and $H$ is a Cartan subgroup of $G$, then manifold $G/H$ is
 the flag manifold ${\mathcal F}$ defined by a parabolic
 subgroup
  can be identified with a
quotient $U_{\mathbb R}/H$, where $U_{\mathbb R}$ is a real
compact form of $G$ and $H$ a closed subgroup of $U_{\mathbb R}$.
  If $\Psi$ is a $1$-dimensional representation of $H$,   then
 $\#\,\Psi(Z(U_{\mathbb R}))$ is also a lower bound for the cardinal of $\pi_1({\mathcal
 H}_{\omega}({\mathcal F}))$ (see {\bf Theorem
 \ref{Thmcartansubgroup}}).

The derivative of the character $\Psi$ is a covector of
${\mathfrak h}$ that can be extended in a trivial way to an
element $\eta\in{\mathfrak u}_{\mathbb R}^*$, using the complement ${\mathfrak l}$ of ${\mathfrak h}$. In the particular
case that the stabilizer of $\eta$ (for the coadjoint action of
$U_{\mathbb R}$) is $H$, the manifold ${\mathcal F}$ is the
coadjoint orbit of $\eta$, the form $\omega$ is the Kirillov
symplectic structure $\varpi$ \cite{Kir0} up to a constant factor
and ${\mathcal H}_{\omega}({\mathcal F})$ is the Hamiltonian group
${\rm Ham}({\mathcal F},\,\varpi)$ \cite{Mc-S}. Under that
hypothesis,
 $\#\,\Psi(Z(U_{\mathbb R}))$ is   a lower bound for the cardinal of $\pi_1({\rm
 Ham}({\mathcal F},\,\varpi))$ ({\bf Corollary \ref{CorHam}}). Here the
 manifold ${\mathcal F}$ is compact and $\omega$ is symplectic, thus, the
 functions $J_Z$ can be normalized. This fact allows us to prove
 the invariance of $\Theta(\psi)$ under deformations of the loop
 $\psi$ in ${\rm Ham}({\mathcal
F},\,\varpi)$; thus, we have a group homomorphism
\begin{equation}\label{Theta:pi_1}
 \Theta:\pi_1({\rm Ham}({\mathcal F},\,\varpi))\longrightarrow
 {\rm U}(1).
  \end{equation}
 From this property, it is possible to obtain an
alternative proof of Corollary \ref{CorHam} (see the last Remark
in Subsection \ref{SubSeccHVF}).

\medskip

\noindent {\bf Case (II).}
 In general,  the Hamiltonian vector fields on quantizable  manifolds \cite{WO} admit lifting to vector fields on
 the prequantum bundles,  which
  are infinitesimal symmetries of the connection.
   We will consider the case of quantizable toric manifolds. Let
   $M$ be such a manifold and $T$ the torus whose action on $M$
   defines the toric structure. Each element ${\bf b}$ in the integer lattice of
   ${\mathfrak t}$ determines a loop
${\varphi_{\bf b}}$ in the
  Hamiltonian group of   $M$. The corresponding gauge transformation $F_1$, on the prequantum bundle,
  can be expressed in a simple way in terms of the mass center  of
  the moment polytope  of the manifold. This fact will allow us to prove {\bf Theorem \ref{ThmCM(Delta)}},
   which gives a sufficient condition for
    ${\varphi_{\bf b}}$  and  ${\varphi_{\bf\tilde b}}$ not be
    homotopic in the Hamiltonian group of $M$. Next, we apply   Theorem \ref{ThmCM(Delta)} to
    the cases when $M$ is the projective space ${\mathbb C}P^n$, a
    Hirzebruch surface and the $1$-point blow up of ${\mathbb C}P^n$.

    If the moment polytope
    of a not necessarily quantizable  toric manifold   has $d$ facets, then the cohomology of the
    manifold is generated by the Chern classes $c_1,\dots,c_d$ of
    $d$ line bundles. When the cohomology class of the symplectic
    structure is of the form $r\sum_j n_jc_j$, with $n_j\in{\mathbb
    Z}$ and $r$ a real number, we will prove {\bf Proposition \ref{Propnonquantiz}}, which is the adaptation  of
    Theorem \ref{ThmCM(Delta)} to this type of
    manifolds.

\smallskip

The paper is organized as follows. In Section \ref{SecP}, from the
representation $\Psi$ of $H$ on a vector space $W$ and a
complement  ${\mathfrak l}$ of ${\mathfrak h}$ in ${\mathfrak g}$, we construct the ${\rm GL}(W)$-principal
fibre bundle ${\mathcal P}$ over $M$ and the $G$-invariant
connection $\Omega$.

An axiomatic definition of the algebras ${\mathfrak X}$ which admit a lift to
${\mathcal P}$ is given in Section \ref{SecX}. In Subsection
\ref{SubSeccLiftcurves}, we study the continuity of the lifting of
certain curves in $M$ to curves in ${\mathcal P}$.
  The lift of isotopies in ${\mathcal G}$ to isotopies in the fibre
bundle, when ${\mathcal G}$
 %being
 is any subgroup of ${\rm Diff}(M)$
whose subalgebra satisfies the axioms required to the subalgebras
${\mathfrak X}$, is studied in Subsection
\ref{SubSeccLiftIsotopies}. In Subsection \ref{SubSeccSubg}
  we prove
 %Propositions \ref{fconstant}, \ref{PropH10} and
 Theorem \ref{pi1G}.

Section \ref{SecHamilt} is
%Subsection \ref{SubSectdim=1}
concerned only with the   case when $\Psi$ is $1$-dimensional.
 In Subsection
\ref{SubSeccAction}, we give an explicit expression for
$\Theta(\psi)$ and discuss its invariance under deformations of
the loop $\psi$.
%the starting point is a given group ${\mathcal G}$
%such that ${\rm Lie}({\mathcal G})$ admits a lift to ${\mathcal
%P}$.
 In Subsection \ref{SubSeccHVF},   we introduce the
algebra ${\mathfrak X}_{\omega}$
 %which admits a lift to ${\mathcal P}$
 and  prove Theorems
\ref{PropHomega} and \ref{Thmcartansubgroup}.

 Section \ref{SectToric} treats with  Case (II). We consider subgroups ${\mathcal G}$ of
${\rm Diff}(M)$, where $M$ is a toric manifold.  In Subsection
\ref{SubsecQuantizable} we assume that the toric manifold is
quantizable and we prove   Theorem \ref{ThmCM(Delta)}. Subsection
\ref{SubsecNonQuantizable} concerns with toric
    manifolds not necessarily quantizable. In this subsection we prove   Proposition \ref{Propnonquantiz}.

%%%%%%%%%%%%%%%%%%%%%%%%%%%%%%%%%%%%%%%%%%%%%%%%%%%%%%%%%%%%%%%%%%%%%%%%%%%%%%%%%%%%%%%%%%%%%%%%%%%%%%%%%

\medskip
\noindent
 {\bf Other approaches.} The results presented in this
article are a small contribution to the study of the homotopy of
subgroups of the diffeomorphism group   of the homogeneous spaces and the  toric manifolds,
   specially the homotopy of   Hamiltonian groups.

 The
homotopy type of the Hamiltonian groups of symplectic manifolds is
only completely known in a few particular cases: When the
dimension of the manifold is $2$, and when the manifold     is a
ruled complex surface \cite{McD1}. In the first case,
 the connected  component of
the identity map in the diffeomorphism group     is homotopy
equivalent to the component of the identity in the
symplectomorphism group;   using this fact,  the topology of the
Hamiltonian group of surfaces can be deduced from the one of the
diffeomorphism group  (see \cite[page 52]{lP01}). In the second case,
 the positivity of   intersections of $J$-holomorphic spheres
in 4-manifolds played a crucial role in the proof of results about
the homotopy type of the corrresponding  Hamiltonian groups  \cite{Ab,A-M,Gr}.

For a general symplectic manifold $N$, given a loop $\psi$ in the
Hamiltonian group of $N$, the Maslov index of the linearized flow
$\psi_{t*}$ gives rise to a map
\begin{equation}\label{MaslovIndex}
\pi_1({\rm Ham}(N))\longrightarrow{\mathbb Z}/2C{\mathbb Z},
\end{equation}
 $C$ being the minimal Chern nunber of $N$
on spheres. Obviously, the map (\ref{MaslovIndex}) can distinguish a maximum of
$2C$ elements in $\pi_1({\rm Ham}(N))$.
 In \cite{V2}, we obtained lower bounds for the cardinal of the
fundamental group of some Hamiltonian groups, by means of the evaluation of the map (\ref{MaslovIndex}) over
certain elements of its domain.
%calculations of certain Maslov indices.
In particular, Theorem 6
of \cite{V2}   gives a lower bound for the cardinal of $\pi_1({\rm
Ham}({\mathcal O}))$, where ${\mathcal O}$ is a coadjoint orbit of
${\rm SU}(n+1)$ diffeomorphic to the quotient of ${\rm
SL}(n+1,\,{\mathbb C})$ by a parabolic subgroup. Of course, that lower bound
is $\leq 2C$.
%, $C$ being the minimal Chern nunber of ${\mathcal O}$ on spheres.
 For the particular case when the parabolic subgroup is
a Borel subgroup, i. e. ${\mathcal O}$ is the flag variety of
$\mathfrak{sl}(n+1,\,{\mathbb C})$,  the mentioned lower bound is
less than or equal to $4$, since the minimal Chern number of this
manifold  is $2$ \cite[page 117]{Mc-S2}. On the other hand,
  the lower bound for $\#\pi_1({\rm
 Ham}({\mathcal F},\,\varpi))$ given in
 Corollary \ref{CorHam} of the present article, when the group $G$ is ${\rm SL}(n+1,\,\mathbb C)$,
depends on $n$ and $\Psi$. In the case that $\Psi$ is injective on
$Z({\rm SU}(n+1))$, this lower bound is $n+1$.

Given a symplectic manifold $N$, there is a group homomorphism
\begin{equation}\label{Idefinition}
I:\pi_1({\rm Ham}(N))\longrightarrow{\mathbb R},
 \end{equation}
  defined through characteristic
classes of the Hamiltonian fibre bundle on  $S^2$ which determine
  each loop in ${\rm Ham}(N)$ (see \cite{V2b}).
 As $I$ is a group
homomorphism, it vanishes on any element of finite order; so, $I$
is not appropriate to detect these elements in  $\pi_1({\rm
Ham}(N))$. But
 studying the homomorphism  (\ref{Idefinition}),  we  proved the existence of infinite
cyclic subgroups in $\pi_1({\rm Ham}(M))$, when $M$ is a toric
manifold (see \cite[Theorem 1.2]{V2b}).
 Also for toric manifolds, McDuff and Tolman proved some results relative to the
corresponding Hamiltonian groups, by developing  the concept of
mass linear pair (see for example
 \cite[Proposition 1.22]{M-T3}).

In \cite{V3}, we considered a linear semisimple Lie group $G$, a
compact Cartan subgroup $T$ and a discrete series representation
$\pi$ of $G$, determined by an element $\phi$ in weight lattice of
${\mathfrak t}$. By means of these data, we constructed a
principal bundle ${\mathcal P}$ on $G/T$ with connection. As here,
 %In the same vein as in the present article
 we considered subgroups of
${\rm Diff}(G/T)$ such that the vector fields in its Lie algebras
admit a lift to vector fields on ${\mathcal P}$, which are
infinitesimal symmetries of the connection.
 We sketched a proof of a result that is a particular case of   Theorem
\ref{pi1G} above mentioned. The present article is, in some
extent, the generalization and formalization of  results sketched
in \cite{V3}.

%%%%%%%%%%%%%%%%%%%%%%%%%%%%%%%%%%%%%%%%%%%%%%%%%%%%%%%%%%%%%%%%%%%%%%%%%%%%%%%%%%%%%%%%%%%%%%%%%%%%%%%%%%%%%%%
\section{The principal bundle ${\mathcal P}$}\label{SecP}

In this section, we introduce some notations that will be used in
the sequel. We will consider pairs $(H,\Psi)$ such
 that, either

 \noindent
 {\it (A)} $\,H$ is an abelian connected closed Lie subgroup of $G$ containing $Z(G)$ and $\Psi$ a
 representation of $H$ in a finite dimensional vector space,

 or

 \noindent
 {\it (B)} $\,H$ is a connected closed Lie subgroup of $G$ containing $Z(G)$ and $\Psi$ a
 $1$-dimensional representation of $H$.

 We denote by
$M$ the quotient $G/H$. Given $A\in{\mathfrak g}$, we put $X_A$
for the vector field on $M$ generated by $A$, through the left
$G$-action.

From now on, we assume that  there exists a vector subspace
${\mathfrak l}\subset {\mathfrak g}$, such that
 \begin{equation}\label{frakl0}
 {\mathfrak
g}={\mathfrak h}\oplus {\mathfrak l}
 \end{equation}
  and $t\cdot
{\mathfrak l}={\mathfrak l}$, for all $t\in H$; here the dot means
the adjoint action. Henceforth, we assume that such a complement
${\mathfrak l}$ of ${\mathfrak h}$ has been fixed. The ${\mathfrak
h}$-component of an element $B\in{\mathfrak g}$ will be denoted
$B_0$. Obviously,
\begin{equation}\label{tB0}
(t\cdot B)_0=t\cdot B_0,\;\hbox{ for all}\; t\in H \; \hbox{and}\;
B\in{\mathfrak g}.
\end{equation}

We will denote by $[{\mathfrak l},\,{\mathfrak l}]_0$ the space
spanned by all elements of the form $[A,\,B]_0$, with $A,B\in
{\mathfrak l}$.

 $\Psi$ is   a group homomorhism $\Psi:H\to{\rm GL}(W)$, where
$W$ is a complex vector space of finite
dimension in the case {\it (A)} and with dimension $1$ in the case
{\it (B)}.  Its derivative will be denoted by $\Psi':{\mathfrak
h}\to\mathfrak{gl}(W)$.

For each $A\in{\mathfrak g}$, we set
 \begin{equation}\label{fagH}
f_A(gH):=\Psi'((g^{-1}\cdot A)_0).
  \end{equation}

  \begin{Prop}\label{PropfA} Formula (\ref{fagH}) defines a map
  $f_A:M\to\mathfrak{gl}(W)$.
  \end{Prop}
  {\it Proof.} In   case {\it (A)}, it is direct consequence of
  (\ref{tB0}) together with the fact $t\cdot B_0=B_0$, for all
  $t\in H$, since $H$ is abelian.

  For the proof in case {\it (B)}, it is sufficient to take into account (\ref{tB0})
  and   that $\Psi'(t\cdot B_0)=\Psi'(B_0)$, since dimension of $\Psi$ is $1$.

 \qed

 %Since $H$ is abelian and
%${\mathfrak l}$ is invariant under adjoint action of $H$, $f_A$ is a
%well-defined map.

 In the case {\it (A)}, from the fact that $\Psi'$ is a Lie
algebra homomorphism and ${\mathfrak h}$ is abelian, it follows that
\begin{equation}\label{[fAfB]}
[f_A,\,f_B]_{\mathfrak{gl}}=0,
 \end{equation}
where $[\,.\,,\,.\,]_{\mathfrak{gl}}$ is the commutator in
$\mathfrak{gl}(W)$. Obviously,   (\ref{[fAfB]}) is also valid  for
the case {\it (B)}.

It is also easy to check that for any $A,C\in{\mathfrak g}$
\begin{equation}\label{XAfC}
 X_A(f_C)=-f_{[A,\,C]}.
\end{equation}

 We denote by
${\mathcal P}$  the following principal ${\rm GL}(W)$-bundle over
$M$
$${\mathcal P}=\big(G\times{\rm GL}(W)\big)/\sim\; \stackrel{\pi}{\longrightarrow} M, $$
where $(g,\,\alpha)\sim(gt,\,\Psi(t^{-1})\alpha)$, for all $g\in
G,$ $\alpha\in{\rm GL}(W)$ and $t\in H$. The projection map from
${\mathcal P}$ to $M$ will be denoted by $\pi$, and the ${\rm
GL}(W)$-right action on ${\mathcal P}$ by $R$.

 $Y_A$ will denote the vector
field   determined by $A\in{\mathfrak g}$ through the left $G$-action on ${\mathcal P}$. The
vertical vector field defined by $y\in\mathfrak{gl}(W)$ will be
denoted $V_y$.

The mapping  $f_A$ can be lifted to a map $ {\bf f}_A$ on ${\mathcal P}$ by
putting
  \begin{equation}\label{bffA}
  {\bf f}_A([g,\,\alpha])=\alpha^{-1}\circ f_A(gH)\circ\alpha,
   \end{equation}
 where $[g,\,\alpha]$ denotes the element of ${\mathcal P}$
  determined by the pair $(g,\,\alpha)\in G\times{\rm GL}(W)$.

   Since $\Psi(t) f_A(gH)=f_A(gH)\Psi(t)$, ${\bf f}_A$ is well-defined.
   From (\ref{bffA}),  it turns out
  \begin{equation}\label{YAbff}
  Y_A({\bf f}_C)=-{\bf f}_{[A,\,C]},\;\;V_y(p)({\bf
  f}_A)=-[y,\,{\bf f}_A(p)]_{\mathfrak{gl}}.
   \end{equation}

In
  case {\it (B)}, (\ref{bffA}) reduces to
   ${\bf f}_A([g,\,\alpha])=f_A(gH)$ and (\ref{YAbff})(ii) to
 $V_y({\bf f}_A)=0$, obviously.

\smallskip

  We denote by $\Omega$ the $G$-invariant connection on ${\mathcal
  P}$ determined by the condition
   \begin{equation}\label{DefConn}
    \Omega(Y_A)={\bf f}_A.
     \end{equation}
  The horizontal lift $X_A^{\sharp}$ of the vector field $X_A$ is
  \begin{equation}\label{XAsharp}
  X_A^{\sharp}([g,\,\alpha])=Y_A([g,\,\alpha])+V_y([g,\,\alpha]),
   \end{equation}
  with  $y:=-{\bf f}_A([g,\,\alpha]).$
  We denote by $D$  the corresponding covariant derivative and by ${\bf K}$ the curvature of this
  connection.  Using the structure equation, (\ref{[fAfB]}), (\ref{YAbff}) and (\ref{DefConn}), it is
  straightforward
to check that
  \begin{equation}\label{bfK}
  {\bf K}(Y_A,\,Y_C)=-{\bf f}_{[A,\,C]}.
  \end{equation}
From (\ref{YAbff}), (\ref{XAsharp}) and (\ref{bfK}), it follows
\begin{equation}\label{DfC}
D{\bf f}_C=-{\bf K}(Y_C,\,.).
 \end{equation}

\medskip

%%%%%%%%%%%%%%%%%%%%%%%%%%%%%%%%%%%%%%%%%%%%%%%%%%%%%%%%%%%%%%%%%%%%%%%%%%%%%%%%%%%%%%%%%%%%%%%%%%%%%%%%%%%%%%%%%%%%
%%%%%%%%%%%%%%%%%%%%%%%%%%%%%%%%%%%%%%%%%%%%%%%%%%%%%%%%%%%%%%%%%%%%%%%%%%%%%%%%%%%%%%%%%%%%%%%%%%%%%%%%%%%%%%%%%%%%%%%%%%%

\section{Algebras admitting a lift to ${\mathcal P}$}\label{SecX}

 From now on,  we assume that the spaces of
 $C^1$ functions between two manifolds are equipped with the Whitney
$C^1$-topology, unless explicit mention is made to the contrary.
In particular, the Lie subgroups of the group of diffeomorphisms
of a manifold and the corresponding Lie algebras will be endowed
with this topology.

\smallskip

\noindent
{\bf Definition 1.} We say that a subalgebra ${\mathfrak X}$  of
${\mathfrak X}(M)$, the algebra consisting of the vector fields on
$M$, admits a lift to ${\mathcal P}$, if   for each
$Z\in{\mathfrak X}$ there is a map ${ a}(Z):{\mathcal
P}\to\mathfrak{gl}(W)$, satisfying the following conditions:

\noindent
 (i) ${ a}(Z)(p\beta)=\beta^{-1}\circ { a}(Z)(p)\circ\beta$, for all
 $\beta\in{\rm GL}(W)$ and all $p\in{\mathcal P}$.

 \noindent
 (ii) $D{ a}(Z)=-{\bf K}(Z^{\sharp},\,.\,),$ where $Z^{\sharp}$ is
 the horizontal lift of $Z$.

  \noindent
 %(iii) For any open set $O$ of $M$, the map $Z_{|O}\mapsto { a}(Z)_{|\pi^{-1}(O)}$ defined on
 %$\{Z_{|O}\,|\, Z\in\mathfrak X\}$ is continuous.
 (iii) The map
 $$Z\in{\mathfrak X}\longmapsto a(Z)\in C^{\infty}({\mathcal
 P},\,\mathfrak{gl}(W))$$
 is continuous with respect to the $C^0$ Whitney topologies in ${\mathfrak X}$ and in $C^{\infty}({\mathcal
 P},\,\mathfrak{gl}(W))$.

\noindent
 (iv) If $X_A\in{\mathfrak X}$, then ${  a}(X_A)={\bf f}_A$.

\smallskip

Note that,  by (\ref{XAsharp}) and (\ref{DfC}),  condition (ii) is
consistent with (iv). Moreover, the Lie algebra
\begin{equation}\label{XA}
\{X_A\,|\,A\in{\mathfrak g}\}
\end{equation}
admits a lift to ${\mathcal P}$.

\smallskip

Given $Z\in{\mathfrak X}$, we set $U(Z)$ for the vector field
 \begin{equation}\label{U(Z)}
U(Z):=Z^{\sharp}+V_{ a(Z)}.
 \end{equation}

 The following proposition is a consequence of condition (i) in
  the above  definition.
  \begin{Prop}\label{Proprighttrans}
   The vector field $U(Z)$ is invariant under the
  right translations in ${\mathcal P}$; that is
  $(R_{\beta})_*(U(Z))=U(Z),$ for all $\beta\in{\rm GL}(W)$.
 \end{Prop}

Since $\Omega$ is a pseudo tensorial form of type ${\rm Ad}$, from
Proposition \ref{Proprighttrans}, it follows
\begin{equation}\label{OmegaPseudo}
 \Omega(U(Z))={\rm Ad}(\beta^{-1})(\Omega(U(Z))),
\end{equation}
 for any $\beta\in{\rm GL}(W)$.

%It is not hard to prove the proposition:
 \begin{Prop}\label{ProprighttransOmega}
The Lie derivative of $\Omega$ along $U(Z)$ vanishes.
 \end{Prop}

{\it Proof.} We put $U$ for $U(Z)$ and $a$ for $a(Z)$. By Cartan's
formula, it is sufficient to check that
  $${\rm
d}(\iota_U\Omega)+\iota_U({\rm d}\Omega)=0.$$
 By (\ref{U(Z)}), ${\rm d}(\iota_U\Omega)={\rm d}a.$ Thus, we need
 to prove that $1$-form on ${\mathcal P}$
 \begin{equation}\label{daiota}
 {\rm d}a+\iota_U({\rm d}\Omega)
  \end{equation}
 vanishes.

 Given an arbitrary
point $p\in{\mathcal P}$,
 we will prove that the $1$-form (\ref{daiota}) applied to a
 horizontal vector $E$ at $p$ vanishes. By (ii) in Definition 1, ${\rm d}a(E)=-{\bf
 K}(Z^{\sharp},\,E)$. Next, we extend $E$ to a horizontal field;
 by the structure equation
 $$\big(\iota_U({\rm d}\Omega)  \big)(E)={\bf K}(U,\,E)-[\Omega(U),\,\Omega(E)]={\bf
 K}(Z^{\sharp},\,E).$$
 So, (\ref{daiota}) vanishes on any horizontal vector.

 Now consider the vertical vector $V_y(p)$, with
 $y\in\mathfrak{gl}(W)$. Again by the structure equation together
 with (\ref{U(Z)})
$$(\iota_U({\rm d}\Omega))(V_y(p))=-[\Omega_p(U),\,\Omega(V_y(p))]=
 -[a(p),\,y].$$
 On the other hand, by (i) in Definition 1,  ${\rm d}a(V_y(p))=-[y,\,a(p)]$. Thus,
 (\ref{daiota}) also vanishes on vertical vectors.
 \qed

%%%%%%%%%%%%%%%%%%%%%%%%%%%%%%%%%%%%%%%%%%%%%%%%%%%%%%%%%%%%%%%%%%%%%%%%%%%%%%%%%%%%%%%%%%%%%%%%%%%%%%%%%%%%%%%%%%%%%%%%%%

\subsection{Lifting of curves in $M$ to curves in ${\mathcal
P}$.}\label{SubSeccLiftcurves}

All the curves on $M$ and on ${\mathcal P}$ considered in this
subsection are defined on the interval $[0,\,1]$. Furthermore,
${\mathfrak X}$ denotes a given algebra that admits a lift to
 ${\mathcal P}$.

\smallskip

 Let $\xi(t)$ be a  curve in $M$ with $\xi(0)=\bar e:=eH$ and $p$ a
 given point in $\pi^{-1}(\bar e)$. Let us assume that there exists a vector field
 $Z\in{\mathfrak X}$,
 such that $Z(\xi(t))=\Dot\xi(t)$, for all $t$.
 Then we can consider the integral curve $\tilde \xi(t)$ of $U(Z)$
 such que $\tilde\xi(0)=p$. By (\ref{U(Z)}), the curve $\tilde\xi(t)$ in ${\mathcal P}$
 can be considered as a lifting (no horizontal lift) of $\xi(t)$.
 In the following paragraphs, we will show the continuity of this
 lifting as a consequence of  property (iii) in
 Definition 1.

\smallskip

%%%%%%%%%%%%%%%%%%%%%%%%%%%%%%%%%%%%%%%%%%%%%%%%%%%%%%%%%%%%%%%%%%%%%%%%%%%%%%%%%%%%%%%%%%
%%%%%%%%%%%%%%%%%%%%%%%%%%%%%%%%%%%%%%%%%%%%%%%%%%%%%%%%%%%%%%%%%%%%%%%%%%%%%%%%%%%%%%%%%%%%%%
Let $\{\xi^c\}_{c\in [0,\,1]}$ be a deformation of $\xi=\xi^0$.
That is, a family of curves  $\xi^c:[0,\,1]\to M$ with initial
point at $\bar e$ (i.e., $\xi^c(0)=\bar e$), satisfying the
following conditions:

(a) The vectors   $\Dot\xi^c(t)=\frac{{\rm d}}{{\rm d}t}\xi^c(t)$
are elements of a vector field $Z^c\in{\mathfrak X}$; that is,
$Z^c(\xi^c(t))=\Dot\xi^c(t)$, for all $t$.

(b) $\{\xi^c\}_{c\in[0,\,1]}$ is a continuous family in the set
$C^1([0,\,1],\,M)$ endowed with the $C^1$ Whitney topology.

For the sake of simplicity, we assume that there exists a chart
$(V,\,\rho)$ on $M$, such that $\xi^c([0,\,1])\subset V$, for all
$c\in[0,\,1]$. The condition (b) involves:

($\alpha$) For any $\epsilon_1>0$, there exists $\delta_1>0$ such
that for all $t\in[0,\,1]$
$$||\xi^c(t)-\xi(t)||<\epsilon_1,\;\;{\rm if}\;\, c<\delta_1.$$
Here $\xi^c(t)$ denotes the coordinates of the point $\xi^c(t)$ in
the chart $(V,\rho)$.

($\beta$) For any $\epsilon_2>0$, there exists $\delta_2>0$ such
that for all $t\in[0,\,1]$
$$||Z^c(\xi^c(t))-Z(\xi(t))||<\epsilon_2,\;\;{\rm if}\;\,c<\delta_2.$$
Here $Z^c(\xi^c(t))$ means the coordinates of the corresponding
vector in the mentioned chart.

From property (a), together with the fact that the curves are
defined on a compact interval and that $\xi^c(0)=\bar e$, one
deduces the equivalence  between properties ($\alpha$) and
($\beta$).

%Since the curves are defined on a compact interval, together with
%property (a) and the fact that $\xi^c(0)=\bar e$, conditions
%($\alpha$) and ($\beta$) are equivalent.

By $\tilde\xi^c$ we denote the lift of $\xi^c$ at the point $p$.
That is, the integral curve of $U^c:=U(Z^c)$ with initial point at
$p$.

The continuity of   the {\it horizontal} lift of  vectors fields
on $M$ to vector fields in ${\mathcal P}$, together with the
condition (iii) in Definition 1 and ($\beta$), imply the following
property: For any $\epsilon
>0$, there exists $\delta >0$ such that for all $t\in[0,\,1]$
$$||U^c(\tilde\xi^c(t))-U(\tilde\xi(t))||<\epsilon,\;\;{\rm if}\;\,c<\delta.$$
 In this formula $U^c$ and $U$ denote the coordinates of the
 respective vector fields in a chart of ${\mathcal P}$ that
 contains the curves $\{\tilde\xi^c(t)\}$.

On the other hand, an equivalence similar to the one between
($\alpha$) and ($\beta$) holds for the curves $\tilde\xi^c$ in
${\mathcal P}$ and the vectors fields $U^c$.

In summary, we have proved the following proposition:
 \begin{Prop}\label{ContLifting}
  With the above notations,
%  fixed $p\in\pi^{-1}(\bar e)$ and
 given a neighborhood  ${\mathcal O}$ of $\tilde \xi$ in the $C^1$ Whitney
 topology of   $C^{1}([0,\,1],\,{\mathcal P})$,
 there is $\delta>0$
such that $\tilde\xi^c\in {\mathcal O}$,
%$$\{\tilde\xi^c(t)\,|\,t\in[0,\,1]\},$$
 %are contained in $O$,
 for all $c<\delta.$
 \end{Prop}

\medskip

%%%%%%%%%%%%%%%%%%%%%%%%%%%%%%%%%%%%%%%%%%%%%%%%%%%%%%%%%%%%%%%%%%%%%%%%%%%%%%%%%%%%%%%%%%%%%%%%%%%%%%%%%%%%%%%%%%%%%%%

\subsection{Lift of isotopies in $M$ to isotopies in ${\mathcal
P}$.}\label{SubSeccLiftIsotopies}
 Let $\{Z_t\}_{t\in[0,\,1]}$ be a
time-dependent vector field on $M$, with $Z_t\in\mathfrak{X}$. Let
$\psi_t$ be the isotopy of $M$ defined by
\begin{equation}\label{isotopypsi}
\frac{{\rm d}\psi_t}{{\rm d}t}=Z_t\circ\psi_t,\;\;\;\;\psi_0={\rm
Id}_M.
\end{equation}

We have the time-dependent vector field $U_t:=U(Z_t)$ on ${\mathcal
P}$ and the corresponding flow ${ F}_t$:
 %$H_t$:
\begin{equation}\label{bdH}
 \frac{{\rm d}{ F}_t(p)}{{\rm d}t}=U_t({ F}_t(p)),\;\;\; {
 F}_0={\rm Id}_{\mathcal P}.
 \end{equation}

By the definition of $U(Z_t)$, it follows that $\pi({
  F}_t(p))=\psi_t(\pi(p))$, for all $t$, where  $p$ is any point of
  ${\mathcal P}$. Thus, if we put $\xi(t)=\psi_t(\bar e)$, its lift
  %$\tilde\xi_t$ of $\xi_t$
   at the point $p\in\pi^{-1}(\bar e)$ is the curve $\tilde\xi(t)=F_t(p)$.

From  Proposition \ref{ProprighttransOmega}, we deduce that
 the diffeomorphism
$F_t$
 preserves the connection form    $\Omega$;    that is, $F_t^*\Omega=\Omega.$

\medskip

   Let us assume that the above isotopy $\psi=\{\psi_t\}_{t\in[0,1]}$ satisfies  $\psi_0=\psi_1={\rm Id}_M$; that
   is,
 $\psi$ is a loop of diffeomorphisms of $M$.
 Then
$F_1$ is a diffeomorphism of ${\mathcal P}$ over the identity,
i.e.,  a gauge transformation of ${\mathcal P}$.
 Thus, there exists a map $h:{\mathcal P}\to{\rm GL}(W)$ satisfying
 $h(p\beta)=\beta^{-1}h(p)\beta$ and  such that
 \begin{equation}\label{bHtf}
 { F}_1(p)=ph(p),\;\; \hbox{for all}\; p\in{\mathcal P}.
 \end{equation}
 Moreover, $F_1^*(\Omega)=h^{-1}{\rm d}h+{\rm Ad}(h^{-1})\circ
\Omega.$

 %By Proposition \ref{ProprighttransOmega},
 % \ref{HtOmega},
Since $\Omega={ F}_1^*(\Omega)$,
 % Hence,
 \begin{equation}\label{bfH1Omeg}
  \Omega =h^{-1}{\rm d}h+{\rm Ad}(h^{-1})\circ
\Omega.
 \end{equation}
%On the other hand, from Proposition
 % \ref{Uinvariante},
 %\ref{Proprighttrans}, it follows
%$\Omega(U_t)={\rm Ad}_{\beta^{-1}}(\Omega(U_t))$, for all
%$\beta\in{\rm GL}(W)$.
 Equality (\ref{bfH1Omeg}) applied to $U_t$ together with
(\ref{OmegaPseudo}) give
%$${\rm Ad}_{\beta^{-1}}(\Omega(U_t))= h^{-1}{\rm d}h(U_t)+{\rm Ad}(h^{-1})(\Omega(U_t)).$$
 %Given $p\in{\mathcal P}$, taking $\beta=h(p)$, we conclude that
 \begin{equation}\label{dhUt}
 {\rm d}h(U_t)=0.
  \end{equation}
  Thus, we have proved the following proposition:
  \begin{Prop}\label{fconstant}
   Let $\psi=\{\psi_t\}_{t\in[0,\,1]}$ be  a loop of diffeomorphisms of
  $M$ generated by a time-dependent vector field $Z_t\in{\mathfrak
  X}$. The map $h:{\mathcal P}\to{\rm GL}(W)$ defined by
  (\ref{bHtf}) is constant on the orbit
  $\{F_t(p)\,|\,t\in[0,\,1]\}$ of any point $p\in{\mathcal
  P}.$
 \end{Prop}

 \smallskip

Let $\{g_t\,|\,t\in[0,\,1]\}$ be a $C^1$ curve   in $G$ with
initial point at the identity element $e$; for brevity such a
curve will be called a {\em path} in $G$. It determines the
isotopy $\varphi=\{\varphi_t\}_{t\in[0,\,1]}$ of $M$, with
 \begin{equation}\label{varphi-t}
\varphi_t(gH)=g_tgH.
\end{equation}
 We put $\{A_t\}$ for the corresponding velocity curve; that is,
$A_t\in{\mathfrak g}$ is defined by
\begin{equation}\label{velocity}
A_t=\Dot g_tg_t^{-1}.
\end{equation}

 Let us assume that $X_{A_t}\in{\mathfrak X}$. This family
gives rise to the time-dependent vector field on ${\mathcal P}$,
$U_t:=U(X_{A_t})$,
 %=Y_{A_t}$,
 which in turn defines a flow ${\bf F}_t$   on ${\mathcal P}$ (the flows in ${\mathcal P}$ determined by paths in $G$ will be denoted in
 boldface). Obviously, the lift of the curve $\varphi_t(\bar e)$ at the point
$p\in\pi^{-1}(\bar e)$ is ${\bf F}_t(p).$
  From (iv) in Definition 1 together with (\ref{XAsharp}), it
  follows
 % (\ref{HorizLift}),
  $U(X_{A_t})=Y_{A_t}$. So,
\begin{equation}\label{FlowH0}
 \frac{{\rm d}{\bf F}_t}{{\rm d}t}=Y_{A_t}\circ {\bf F}_t,\;\;\;
 {\bf F}_0={\rm Id}_{\mathcal P}.
 \end{equation}
 \begin{Lem}\label{LemHt}
 The bundle diffeomorphism ${\bf F}_t$ defined in (\ref{FlowH0}) is
the left multiplication by $g_t$ in ${\mathcal P}$.
 %; i.e. $H_t={\mathcal L}_{g_t}.$
\end{Lem}

{\it Proof.} Given $p\in{\mathcal P}$, by (\ref{velocity})
 $$ \frac{{\rm d}}{{\rm
d}u}\bigg|_{u=t}g_up= \frac{{\rm d}}{{\rm
d}u}\bigg|_{u=t}g_ug_t^{-1}g_tp=Y_{A_t}(g_tp).$$
 \qed

\smallskip

 If $g_1$ is an element of the center of $G$, then $\varphi_0=\varphi_1$ and
 ${\bf F}_1$
 %, the time-$1$ map of that flow,
  is a gauge transformation of ${\mathcal P}$.

\smallskip

Henceforth, we assume that {\em $\Psi(g_1)$ is a scalar operator
for any $g_1\in Z(G)$.} This condition holds trivially in case
{\it (B)}.   The case {\it (A)} with this additional condition
will be designed as {\it (A')}.

\begin{Prop}\label{PropH10}
Let $g_t$ be a path in $G$, such that $X_{A_t}\in{\mathfrak X}$,
where $A_t=\Dot g_tg_t^{-1}$. If  $g_1\in Z(G)$, then ${\bf
F}_1(p)=p\Psi(g_1)$,     for all $p\in{\mathcal P}$.
 \end{Prop}
 {\it Proof.} By Lemma \ref{LemHt},
$${\bf F}_1([g,\,\alpha])=[g_1g,\,\alpha]=[gg_1,\,\alpha]=[g,\,\Psi(g_1)\alpha]=[g,\,\alpha]\Psi(g_1),$$
since $\Psi(g_1)$ is a scalar operator.
 \qed

\smallskip

{\bf Remark.}    Note that  paths in $G$ with the same final point
in $Z(G)$ determine the same gauge transformation ${\bf F}_1$.

\smallskip

%%%%%%%%%%%%%%%%%%%%%%%%%%%%%%%%%%%%%%%%%%%%%%%%%%%%%%%%%%%%%%%%%%%%%%%%%%%%%%%%%%%%%%%%%%%%%%%%%%%%%%%%%%%%%%

\subsection{Subgroups of ${\rm Diff}(M)$.}\label{SubSeccSubg}

Let ${\mathcal G}$ be a connected Lie subgroup of
 ${\rm Diff}(M)$  such that:
\begin{equation}
\begin{aligned}\label{aligned}
 &\hbox{(i)} \;{\rm Lie}( {\mathcal G}) \;\hbox{admits a lift to}\;{\mathcal P}. \\
 & \hbox{(ii)} \;\hbox{If} \; \{g_t\} \; \hbox{is a path in}\; G,\;
 \hbox{then the isotopy}\;
\{\varphi_t\} \;
%\hbox{defined in}\;(\ref{varphit})\;
 % \hbox{with}\; \varphi_t(gT_{\mathbb R})=g_tgT_{\mathbb R}\; \hbox{is an isotopy that  belongs
 \hbox{is contained in}\;
   {\mathcal G}.
\end{aligned}
\end{equation}
 We will  prove that
$$\#\,\pi_1({\mathcal G})\geq\#\,\{\Psi(g)\,|\,g\in Z(G)
\},$$
 but we need to introduce some notations.

Let $\{g_t\}$ be a path in  $G$ such that
  $g_1\in Z(G)$.  Let
$\{\zeta^s\}_s$ be a {\em deformation} of the loop $\varphi$ in
${\mathcal G}$. That is, for each $s$, $\zeta^s$ is a loop in
${\mathcal G}$ at ${\rm Id}_M$, with $\zeta^0=\varphi$. We also
assume that $s\mapsto \zeta^s$ is a continuous map (considering
${\mathcal G}$ equipped with the $C^1$-topology, as we said).
%, the topology of uniform convergence of all
%derivatives.

For each $s$, we have the corresponding time-dependent vector
field $Z^s_t\in{\rm Lie}({\mathcal G})$, given by
$$\frac{{\rm d}\zeta^s_t}{{\rm d}t}=Z^s_t\circ\zeta^s_t.$$
The respective time-dependent vector field on ${\mathcal P}$,
$U^s_t:=(Z^s_t)^{\sharp}+V_{a(Z^s_t)}$ determines the
corresponding flow $F^s_t$. As above, $F^s_1$ is a gauge
transformation, which can be written
\begin{equation}\label{Fs1}
 F^s_1(q)=q h^s(q), \; \hbox{for all} \; q\in{\mathcal P},
 \end{equation}
 with $h^s$ constant  on the orbits
 $\{F^s_t(q)\,|\,t\in[0,\,1]\}$.
%By Proposition \ref{fconstant},
% $h^s$ is constant on
%\begin{equation}\label{Fst}
%    \{{ F}^s_t(q)\,|\, t\in[0,\,1]\},
%     \end{equation}
%the orbit of any $q$. That is, there exists an element
%$\beta\in{\rm GL}(W)$ such that
%\begin{equation}\label{hs}
%h^s(p')=p'\beta,\;\;\;\hbox{for all}\;p'\in\{F^s_t(q)\,|\,
%t\in[0,\,1]\}.
%\end{equation}

 Fixed $s$, we can also consider the loop $\xi^s$ in $M$, obtained by evaluating $\zeta^s_t$ at the point $\bar e$
\begin{equation}\label{xist}
 \{\xi^s(t):=\zeta^s_t(\bar e)\,|\, t\in[0,\,1]  \}.
\end{equation}
  Obviously, the lifting  of this loop to ${\mathcal P}$ at a point $p\in\pi^{-1}(\bar e)$
  is just the curve $F^s_t(p)$; i.e.,
\begin{equation}\label{tildexis}
\tilde\xi^s(t)=F^s_t(p).
\end{equation}

\smallskip

 In the statement of the following theorem we refer to the $H$-principal bundle
 \begin{equation}\label{TRbundle}
H\to G\stackrel{{\rm pr}}{\longrightarrow} M=G/H.
 \end{equation}
  We assume that this bundle is endowed  with the invariant connection
 \cite[page 103]{Kob-Nom}
 determined by the splitting (\ref{frakl0}).
 % ${\mathfrak g}={\mathfrak h}\oplus{\mathfrak l}$.
 It is well-known that
 the Lie algebra  of the holonomy  group ${\rm
Hol}_e$ at $e\in G$ (of this invariant connection)  is $[\mathfrak
l,\,\mathfrak l]_0$
%generated by the vectores of the form
%$[A,\,C]_0$, with $A,C\in{\mathfrak l} $
 (see \cite[Theorem 11.1,
page 103]{Kob-Nom}).

 \begin{Thm}\label{ThmPerturb}
Let $\{g_t\}_{t\in[0,\,1]}$ be a path in $G$ which
 is a horizontal curve with
respect the invariant connection in $G \stackrel{{\rm
pr}}\longrightarrow M$ and such that its final point  belongs to
$Z(G)$.  If $\{\zeta^s\}_s$ is a deformation of the loop
 %$\{\varphi_t\}$
 $\varphi$ in ${\mathcal
G}$, then the gauge transformation ${F}^s_1$ of ${\mathcal P}$
defined by $\zeta^s$ satisfies
 $${ F}^s_1(p)=p\Psi(g_1),$$
 for all $p\in\pi^{-1}(\bar e).$
 %\;\; \hbox{for all}\;\; p\in\pi^{-1}(\bar e).$$
% with $\kappa(g_1)$ defined in (\ref{Psig1}).
  %$\kappa(g_1)\,{\rm Id}_{W}=\Psi(g_1)$.
 \end{Thm}

{\it Proof.} Fix an arbitrary point $p\in{\mathcal P}$, such that
$\pi(p)=\bar e$. The idea of the proof is to construct a path in $G$, with
final point at $g_1$, that defines
 a loop in $M$ as close   to the evaluation curve $\xi^s$ (see (\ref{xist})) as we wish. Then the lifting of the loop
 at the point $p$  and the one of $\xi^s$,  give rise to the constants $\Psi(g_1)$ and
$h^s(p)$, whose difference will be as small as we wish.
 % which  will be ``as close as we wish".

 We have the following closed curves in
$M$:
 \begin{equation}\label{curves}
 \{\varphi_t(\bar e)\,|\,t\in[0,\,1]\},\;\;\;
 \xi^s\equiv\{\xi^s(t) \,|\,t\in[0,\,1]\}.
  \end{equation}

\smallskip

 {\it (I) Horizontal lifting from $M$ to $G$.}
 The invariant connection on the principal bundle $G\to
 M$ determines the horizontal lifting of the curves (\ref{curves}).
 We denote by  $\upsilon^a$  the curve in $G$  horizontal lift
 of
 $\xi^a$
 %$\{\zeta^a_t(\bar e)\,|\, t\in[0,\,1]\}$
 at the point $e\in G$.
 By hypothesis, $g_t$ is the   horizontal lift
  of $\{\varphi_t(\bar e)=g_t\bar e\}$  at the point  $e$. In particular,
  $\upsilon_t^0=g_t$.

 Fixed $s$, a value of the parameter of deformation, since the horizontal lifting is a continuous
 operation,
 $$S:=\{\upsilon^a_t\,|\,a\in[0,\,s],\;t\in[0,\,1]\}$$
 is a piece of surface in   $G$. The element $g_1$, final point of $\upsilon^0$, belongs to this
connected surface in $G$.

 Given $0<\epsilon<1,$
let $\{\mu_t\,|\,t\in[1-\epsilon,\,1]\}$ be a curve in $S$ with
$\mu_{1-\epsilon}=\upsilon^s_{1-\epsilon}$ and final point  at
$g_1$. The curve $\mu$ can be chosen so that its juxtaposition
with the restriction of $\upsilon^s$ to $[0,\,1-\epsilon]$ is a ${
C}^1$ curve in $G$ (see Figure 1). The path in $G$ defined by this
juxtaposition will be denoted  $\{
  g'_t\,|\, t\in[0,\,1]\}$; its final point is $g_1$.

By taking $\epsilon$ sufficiently small we can get a curve $\{
g'_t\bar e={\rm pr}( g'_t)\,|\,t\in[0,\,1]\}$ in $M$ contained in
an arbitrary   neighborhood of the curve $\xi^s$ in the $C^1$
topology of $C^1([0,\,1],\,M)$.
%$\{\xi^s(t)\,|\,t\in[0,\,1]\}$.
%$\{\zeta^s_t(\bar e)\,|\, t\in[0,\,1]\}$ (see Figure 2).
 Thus, we have proved the following lemma:
\begin{Lem}\label{LemMainTh}
Given a    neighborhood of the curve $\xi^s$ in the $C^1$ topology
of $C^1([0,\,1],\,M)$,
% $\{\xi^s(t)\,|\,t\in[0,\,1]\}$,
there is a path $\{g'_t\}$ in $G$, with final point at $g_1$, such
that the loop $\{g'_t\bar e\}$ is contained in that neighborhood.
 \end{Lem}

\begin{figure}[htbp]
\begin{center}
\epsfig{file=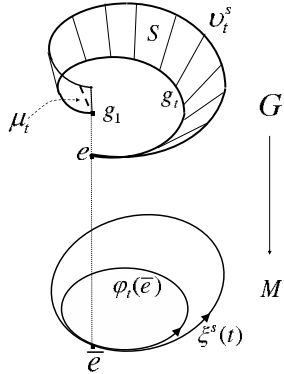, height=5.5cm}
 \caption[Figure 1]
 {\small The path $\{g'_t\}$ is the juxtaposition of   $\upsilon^s_{\,|[0,\,1-\epsilon]}$ and   $\mu_t$.}
 %\Delta_2$ bundle over $\Delta_1$}
\end{center}
\end{figure}

\smallskip

{\it (II) Lift from $M$ to ${\mathcal P}$.}
 As we said,  $p$ will be an arbitrary point of $\pi^{-1}(\bar e)$. The loop of diffeomorphisms $\zeta^s$ defines the
    corresponding flow $F^s_t$
   and $F^s_1$ is determined by   $h^s$ (see
   (\ref{Fs1})).

We can also consider $\tilde \xi^s$, the lifting  to ${\mathcal
P}$ at the point $p$ of the loop in  $\xi^s$.
 %defined in (\ref{xist});
 By (\ref{Fs1}) and (\ref{tildexis})
 \begin{equation}
 \tilde\xi^s(1)=F^s_1(p)=ph^s(p).
 \end{equation}

%$$\xi^s_t=\{\zeta^s_t(\bar e)\,|\, t\in[0,\,1]  \}.$$
%Obviously, this lifting is just the curve (\ref{Fst}); i.e.,
%$\tilde\xi^s_t=F^s_t(p)$.

   % This curve  the lift to ${\mathcal P}$ at the point $p$ of the loop
    %$\{\zeta^s_t(\bar e)\,|\, t\in[0,\,1]  \}$ in $M$.

    On the other hand, the path $\{g'_t\}$ has final point  at
    $g_1\in Z(G)$. It determines  the respective flow ${\bf F}'_t$ and, according to Proposition
    \ref{PropH10},
    %\ref{PropH1},
 ${\bf F}'_1$ is the gauge transformation defined by the constant
 $\Psi(g_1)$
 \begin{equation}\label{h'}
{\bf F}'_1(q)=q\Psi(g_1).
  \end{equation}
   %(Proposition \ref{Psi(g1)=kappa}).
We have also the lifting to ${\mathcal P}$ at the point $p$ of the
loop $\{g'_t\bar e \}$; this lift is, of course, the curve
 \begin{equation}\label{HatHt}
\{{\bf F}'_t(p)\,|\,t\in[0,\,1]   \}
  \end{equation}

\smallskip

{\it (III) Continuity of the lifting.}
 Let ${\mathcal O}$ be
an arbitrary   neighborhood of the curve $\tilde\xi^s$ in the
$C^1$ topology of $C^1([0,\,1]\,{\mathcal P})$. By the continuity
of the lifting stated in Proposition \ref{ContLifting}, together
with Lemma \ref{LemMainTh},
 there are curves of the form $\{g'_t\bar e\}$,
 where $\{g'_t\}$ is a path in $G$ with final point at $g_1$, whose
lifting (\ref{HatHt}) to
%${\mathcal P}$
belong to ${\mathcal O}$. Hence, $\tilde\xi^s(1)=ph^s(p)$ is a
point of ${\mathcal P}$ as close to one point of the type ${\bf
F}'_1(p)$ as we wish. But, by (\ref{h'}), all the points ${\bf
F}'_1(p)$ are equal to $p\Psi(g_1)$ (independently of the path
$\{g'_t\}$).  Therefore, the constant $h^s(p)$ is $\Psi(g_1)$;
that is, $F^s_1(p)=p\Psi(g_1)$, for all $p\in\pi^{-1}(\bar e)$.
 \qed

\begin{Cor}\label{CorPert}
 Let $\{g_t\}$ and $\{\tilde g_t\}$ be
 paths in $G$ satisfying the hypotheses of Theorem
 \ref{ThmPerturb}. If $\Psi(g_1)\ne\Psi(\tilde g_1)$, then the
 corresponding loops $\varphi$  and $\tilde\varphi$   are not
 homotopic in ${\mathcal G}$.
 \end{Cor}

If  $Z(G)\subset {\rm Hol}_e$,  each element of $Z(G)$ can be
joined to $e$ by a horizontal curve in $G$ (horizontal with
respect to the invariant connection).

The following proposition gives a sufficient condition for
$Z(G)\subset{\rm Hol}_e$, when ${\mathfrak g}$ is a perfect Lie
algebra; i.e. ${\mathfrak g}=[{\mathfrak g},\,{\mathfrak g}]$.

\begin{Prop}\label{Z(G)Hol}
%Let ${\mathfrak g}$ be a perfect Lie algebra. If $Z(G)\subset H$ and $H$ is abelian,
 %then $Z(G)\subset{\rm Hol}_e$.
  If ${\mathfrak g}$ is a perfect Lie algebra and $H$ is abelian,
 then $Z(G)\subset{\rm Hol}_e$.
\end{Prop}

{\it Proof.} Since ${\mathfrak h}$ is abelian and ${\mathfrak l}$
is invariant under ${\rm Ad}(H)$,
 \begin{equation}\label{[]}
 {\mathfrak h}\oplus{\mathfrak l}={\mathfrak g}=
[{\mathfrak g},\,{\mathfrak g}]={\mathfrak l}_1+[{\mathfrak
l},\,{\mathfrak l}],
 \end{equation}
  with
 ${\mathfrak l}_1\subset {\mathfrak l}$. By taking the ${\mathfrak
 h}$-component in (\ref{[]}), we obtain
 ${\mathfrak h}=[{\mathfrak l},\,{\mathfrak l}]_0$; as $H$ is
 connected, $H\subset{\rm Hol}_e$. By hypothesis,  $Z(G)$ is a subgroup of
 $H$ in case (A) and in case (B).

  \qed

We denote by ${\rm exp}([{\mathfrak l},\,{\mathfrak l}]_0)$ the
connected subgroup of $G$ whose Lie algebra is $[{\mathfrak
l},\,{\mathfrak l}]_0$. {\em Henceforth, we assume that the
following property holds}
\begin{equation}\label{Z(G)subset}
 Z(G)\subset {\rm exp}([{\mathfrak l},\,{\mathfrak l}]_0).
 \end{equation}

From Corollary \ref{CorPert}, it follows the following theorem:
 \begin{Thm}\label{pi1G}
 Let $G$ be a  Lie group, and $(H,\Psi)$ a pair satisfying
  either   condition {\it
 (A')} or
  condition {\it (B)}, and ${\mathfrak l}$ an ${\rm Ad}(H)$-invariant complement of ${\mathfrak h}$
  for which (\ref{Z(G)subset}) is valid.   If ${\mathcal G}$ is a connected Lie subgroup of ${\rm Diff}(M)$ for which
  (\ref{aligned})  holds,
   %and $Z(G)\subset H$,
    then
 $$\#\,\pi_1({\mathcal G})\geq\#\,\Psi(Z(G)).$$
 \end{Thm}

\smallskip

Given $g\in G$, we denote by ${\bf g}$ the diffeomorphism of $M$
defined by left action of $g$; that is, ${\bf g}(bH)=gbH$.  We put
${\bf G}$ for the subgroup of ${\rm Diff}(M)$ consisting of all
diffeomorphisms ${\bf g}$. Obviously,  the group ${\mathcal
G}\subset{\rm Diff}(M)$ satisfies (\ref{aligned})(ii) if ${\bf
G}\subset {\mathcal G}.$

 Each closed path in $G$ gives rise to a loop in
${\bf G}$, and homotopic closed paths in $G$ induce homotopic
loops in ${\bf G}$. However, there are loops in ${\bf G}$ which
are not defined by closed paths in $G$.
 Evidently, ${\rm Lie}({\bf
G})$ is the Lie algebra (\ref{XA}).
%=\{X_A\,|\,A\in{\mathfrak
%g}\}\subset{\mathfrak X}(M)$.
This Lie algebra satisfies the
conditions required to the algebras $\mathfrak X$ corresponding to
any representation $\Psi$ of $H$. Thus, we have the following corollary
of Theorem \ref{pi1G}:
\begin{Cor}\label{CorbfG}
 If $H$, $\Psi$ and ${\mathfrak l}$ satisfy the hypotheses of
 Theorem \ref{pi1G},
 %and $Z(G)\subset H$,
 then
$$\#\,\pi_1({\bf G})\geq \#\,\Psi(Z(G)).$$
 %for any  representation $\Psi$ of $H$ such that the pair $(H,\,\Psi)$ satisfies condition  {\it (B)}.
\end{Cor}

\medskip

%%%%%%%%%%%%%%%%%%%%%%%%%%%%%%%%%%%%%%%%%%%%%%%%%%%%%%%%%%%%%%%%%%%%%%%%%%%%%%%%%%%%%%%%%%
%%%%%%%%%%%%%%%%%%%%%%%%%%%%%%%%%%%%%%%%%%%%%%%%%%%%%%%%%%%%%%%%%%%%%%%%%%%%%%%%%%%%%%%%%%%%%%%%%%%%%%%%%%%%%%%%%%%%%%%%%%
%%%%%%%%%%%%%%%%%%%%%%%%%%%%%%%%%%%%%%%%%%%%%%%%%%%%%%%%%%%%%%%%%%%%%%%%%%%%%%%%%%%%%%%%%%%%%%%%%%%%%%%%%%%%%%%%%%%%%%%%

\section{Hamiltonian groups}\label{SecHamilt}

In this Section {\em we assume that the pair $(H,\,\Psi)$ belongs
to case {\it (B)}}. When $\Psi$ is a $1$-dimensional
representation, the curvature ${\bf K}$ projects a presymplectic
structure $\omega$ on M, as we said in the Introduction.

In Subsection \ref{SubSeccAction}, we will consider the
$\omega$-action integral around the loop of diffeomorphisms
generated by a time-dependent vector field in an algebra which
admits a lift to ${\mathcal P}$. We will relate this action
integral with the gauge transformation associated with that loop.

In Subsection \ref{SubSeccHVF}, we consider the algebra
${\mathfrak X}_{\omega}$ of all $\omega$-Hamiltonian vector fields
on $M$. We will give conditions for ${\mathfrak X}_{\omega}$
admitting  a lift to ${\mathcal P}$ and we will apply Theorem
\ref{pi1G} to this particular case.

 %%%%%%%%%%%%%%%%%%%%%%%%%%%%%%%%%%%%%%%%%%%%%%%%%%%%%%%%%%%%%%%%%%%%%%%%%%%%%%%%%%%%%%%%%%%%%%%%%%%%%%%%%%%%%

\subsection{Action integral.}\label{SubSeccAction}  When ${\rm
dim}\,W=1$,   Proposition \ref{fconstant} can be stated as:
\begin{Prop}\label{fconstant2}
   Let $\psi=\{\psi_t\}_{t\in[0,\,1]}$ be  a loop of diffeomorphisms of
  $M$ generated by a time-dependent vector field $Z_t$ that belongs to a Lie algebra ${\mathfrak
  X}$ admitting a lift to ${\mathcal P}$. If ${\rm
dim}\,W=1$, then there exists a complex number $\Theta(\psi)$ such that,
  \begin{equation}\label{F1(p)=}
  F_1(p)=p\Theta(\psi),
  \end{equation}
  for all $p\in{\mathcal P}.$
 \end{Prop}
{\it Proof.} As dimension of $W$ is $1$, (\ref{bfH1Omeg}) reduces
to $dh=0$. So, $h$ is constant on ${\mathcal P}$.

\qed

 Let $Z$ be a vector field in an algebra
${\mathfrak X}$ that admits a lift to ${\mathcal P}$; now,  by property (i) in Definition 1, $a(Z)$
projects a function $J_Z$ on $M$ defined by
\begin{equation}\label{a(Z):0}
 J_Z(\pi(p))=a(Z)(p).
 \end{equation}
  Analogously, the curvature ${\bf K}$
projects a closed $2$-form $\omega$ on $M$; thus, $\omega$ is a
presymplectic structure on $M$.   The condition (ii) in Definition
1 gives rise to
 \begin{equation}\label{dJZ}
 {\rm d}J_Z=-\omega(Z,\,.\,).
 \end{equation}
 Note that $Z$ determines $J_Z$  by
 (\ref{dJZ}), up to an additive constant. However, if $Z=X_A$, then   condition (iv) in
 Definition 1 implies that $J_{X_A}=f_A$.

Next, we will express the number $\Theta(\psi)$ defined in
Proposition \ref{fconstant2} as   a ``generalized action integral"
on the presymplectic manifold $(M,\,\omega)$.

%As $T$ is compact, we can assume that the representation $\Psi$ is
%unitary.
 Let $\psi$ be the loop considered in Proposition
\ref{fconstant2}. Given $x\in M$, we put $\{x_t\}$ for the closed
evaluation curve $\{\psi_t(x)\}$. For $p\in\pi^{-1}(x)$, we set
$p_t:=F_t(p)$, thus, $\pi(p_t)=x_t$. Let $\sigma$ be a section of
the principal bundle ${\mathcal P}$, such that $x$ belongs to the
domain of $\sigma$ and $\sigma(x)=p$. We define the complex number
$\theta_t$ by the relation
 \begin{equation}\label{mt}
F_t(p)=\sigma(x_t)\theta_t.
 \end{equation}
Obviously $\theta_0=1$ and $\theta_1=\Theta(\psi)$. Deriving
(\ref{mt}) with respect to $t$,
$$\frac{{\rm d}F_t(p)}{{\rm d}t}=(R_{\theta_t})_*\big(\sigma_*(Z_t(x_t))
\big)+V_{\Dot \theta_t/\theta_t}(p_t).$$
  Taking into account that
  $$\sigma_*(Z_t(x_t))=Z^{\sharp}_t(\sigma(x_t))+V_y(\sigma(x_t)),$$
  with $y=\Omega(\sigma_*(Z_t(x_t)))=\alpha(Z_t(x_t)),$ where
  $\alpha$ is the connection form in the trivialization defined by
  $\sigma$. Then,
  \begin{equation}\label{fracdFt1}
   \frac{{\rm d}F_t(p)}{{\rm d}t}=Z^{\sharp}_t(p_t)+V_y(p_t)+V_{\Dot
  \theta/\theta}(p_t).
   \end{equation}
  On the other hand, by (\ref{bdH}),
 \begin{equation}\label{fracdFt2}
  \frac{{\rm d}F_t(p)}{{\rm d}t}=Z^{\sharp}_t(p_t)+V_{a_t}(p_t),
  \end{equation}
  where $a_t:=a(Z_t)(p_t)$.  From (\ref{fracdFt1}),
  (\ref{fracdFt2}) and (\ref{a(Z):0}),
  it follows
   $$\frac{\Dot \theta}{\theta}=-\alpha(Z_t(x_t))+J_{Z_t}(x_t).$$
   Hence,
   $$\Theta(\psi)={\rm exp}\Big(-\int_0^1\alpha(Z_t(x_t))dt+\int_0^1J_{Z_t}(x_t) dt    \Big).$$

If $H_1(M,\,{\mathbb Z})=0$, there exists a $2$-chain $C$ in $M$
whose boundary is the closed curve $x_t$. In this case, by Stokes'
theorem and taking into account that the curvature ${\bf
   K}$ projects on $M$ the closed $2$-form $\omega$, we have
   the following proposition:
   \begin{Prop}\label{PropTheta}
    Let $\psi$ be the loop of Proposition \ref{fconstant2} and $x$ a point of $M$. If
    $H_1(M,\,{\mathbb Z})=0$, then
 \begin{equation}\label{Theta}
    \Theta(\psi)={\rm exp}\Big(-\int_C\omega+\int_0^1 J_{Z_t}(x_t) dt    \Big),
    \end{equation}
where $C$ is a $2$-chain whose boundary is the closed curve
$\{x_t:=\psi_t(x)\}$.
% and  $J_t:=J_{Z_t}(x_t).$
 \end{Prop}
   From the equality (\ref{Theta}), one deduces that the  right hand side of this equation
  is independent of the point $x$. It has the form  of the exponential
of an action integral (see the last Remark of  Subsection \ref{SubSeccHVF}).

  In general, neither $M$ is compact nor $\omega$ is
symplectic. From these negative facts, it turns out that
 the evaluation curve $x_t$ is not necessarily  nullhomologous
\cite[Lemma 10.31]{Mc-S};
 so, it is necessary to impose the
vanishing  of the first homology group of $M$ to write
(\ref{Theta}).

% (b) The impossibility of fixing the  functions
%$J_{Z_t}$ by means a consistent normalization. That is, in the
%axiomatic definition of the subalgebra ${\mathfrak X}$ the
%properties we have impose to   the functions $a(Z)$ do not
%determine uniquely these functions. Thus, the right hand side of
%(\ref{Theta})  depends on  the functions $J_{Z_t}$, which are not
%determined uniquely by $\psi$.

In general $\Theta$ is not invariant under deformations of the
loop $\psi$, essentially because the $J_Z$ are not normalized.
More precisely, let $\kappa^s$ be a family of loops in ${\rm
Diff}(M)$ at ${\rm Id}_M$, each of which satisfies the hypotheses
of Proposition \ref{fconstant2}, and such that $\kappa^0=\psi$.
Then we have the corresponding time-dependent vector fields
$Z^s_t\in {\mathfrak X}$, the maps $J_{Z^s_t}$ with
$J_{Z^0_t}=J_{Z_t}$ and the corresponding    invariants
$\Theta^s:=\Theta(\kappa^s)$ (defined by means of the
$J_{Z^s_t}$)
 $$\Theta^s={\rm exp}\Big(-\int_{C^s}\omega+\int_0^1J_{Z^s_t}(\kappa^s_t(x))\, dt    \Big),$$
where $C^s$ is a $2$-chain such that $\partial C^s$ is the curve
$\{\kappa^s_t(x)\}_t$.

 On the other hand, the variation of $\kappa^s_t(x)$ with
respect to $s$ defines the time-dependent vector field $B_t$
 $$B_t(\kappa^s_t(x)):=\frac{\partial}{\partial
s}\kappa^s_t(x).$$
%Since $\frac{i}{2\pi}\omega$ defines an element of
%$H^2(M,\,{\mathbb Z})$,
Hence,
 $$\frac{1}{\Theta(\psi)}\frac{{\rm d}\Theta^s}{{\rm
d}s}\bigg|_{s=0}=\int_0^1\omega(Z_t(x_t),\,B_t(x_t))\,dt+\int_0^1
B_t(x_t)\big(J_{Z_t}\big)\,dt+\int_0^1\Dot J_t(x_t).$$
 where
\begin{equation}\label{DotJt}
 \Dot J_t:=\frac{{\rm d}}{{\rm d}s}\bigg|_{s=0}J_{Z^s_t}.
 \end{equation}
  Since
$\iota_{Z_t}\omega =-{\rm d}J_{Z_t}$,
\begin{equation}\label{varTheta}
\frac{1}{\Theta(\psi)}\frac{{\rm d}\Theta^s}{{\rm d}s}\bigg|_{s=0}
 =\int_0^1\Dot J_t(x_t)\,dt.
\end{equation}
Thus, we have the following proposition:
\begin{Prop}\label{dThetads}
Let $\{\kappa^s\}_s$ be a deformation of the loop $\psi$, where
each $\kappa^s$ is generated by a time-dependent vector field
belonging to ${\mathfrak X}$. Then
 $$\int_0^1\Dot J_t(x_t)\,dt$$
 is independent of the point point $x$, and it equals
 $$\frac{1}{\Theta(\psi)}\frac{{\rm d}\Theta(\kappa^s)}{{\rm
 d}s}\bigg|_{s=0}.$$
  \end{Prop}

 %As the left hand side of (\ref{varTheta}) is independent of $x$, so is the right hand side.

Let $\varphi$ be   the loop in the statement of Theorem
\ref{ThmPerturb}. After Proposition \ref{fconstant2},  Theorem
\ref{ThmPerturb} asserts (when ${\rm dim}\, W=1$)  the following:
 \begin{Prop}\label{Thmperturbbis}
 Let $\varphi$ be the loop in the statement of Theorem
 \ref{ThmPerturb} and $\{\zeta^s\}_s$ a deformation of $\varphi$ in ${\mathcal
 G}$,  then
 $$\frac{{\rm d}\Theta(\zeta^s)}{{\rm d}s}\bigg|_{s=0}=0.$$
 \end{Prop}

 As we will see in the last Remark at the end of this Section, when
  $G$ is semisimple,
 $M$ is a flag manifold of ${\mathfrak g}$   and $\omega$  the Kirillov  symplectic
 form, the constant (\ref{varTheta}) vanishes for any deformation in the Hamiltonian group of a given
 Hamiltonian loop $\psi$. In more precise terms, in this case, it is possible to modify the maps $a(Z)$, whose existence is postulated in
 the definition of ${\mathfrak X}$, so that: $(\alpha)$ The new maps satisfy
 the conditions imposed in Definition 1.  $({\beta})$ For the new
 maps $J_Z$, integral (\ref{varTheta}) is zero.

\smallskip

%%%%%%%%%%%%%%%%%%%%%%%%%%%%%%%%%%%%%%%%%%%%%%%%%%%%%%%%%%%%%%%%%%%%%%%%%%%%%%%%%

\subsection{Presymplectic structure.}\label{SubSeccHVF}
 By (\ref{bfK}) the $2$-form $\omega$, projection of ${\bf K}$ on $M$, satisfies
 \begin{equation}\label{omega(XA,XB}
\omega(X_A,\,X_B)=-f_{[A,\,B]}.
\end{equation}
 This relation allows to define the form $\omega$ by means the representation $\Psi$  and the complement
  ${\mathfrak l}$ directly, without
 resorting to the principal bundle ${\mathcal P}$.
 From (\ref{omega(XA,XB}), using the Jacobi identity, one can prove directly that
 $\omega$ is a closed form
(see \cite[page 6]{Kir}).

We put
$${\mathfrak X}_{\omega}:=\{Z\in\mathfrak{X}(M)\,|\,
\omega(Z,\,.\,)\;\hbox{is exact}\,\}.$$

Let ${\mathcal H}_{\omega}(M)$ denote the group consisting of the
time-$1$ maps of the  isotopies in $M$ generated by a
time-dependent vector fields $Z_t\in {\mathfrak X}_{\omega}$. If
$\omega$ is a symplectic form, ${\mathcal H}_{\omega}(M)$ is the
corresponding Hamiltonian group ${\rm Ham}(M,\,\omega)$.
% By
%Cartan's formula the Lie derivative of $\omega$ along Z vanishes
%\begin{equation}\label{Lieomega}
%{\mathfrak L}_Z\omega=0
 %\end{equation}

 By (\ref{XAfC}),
 \begin{equation}\label{omega(XA}
\omega(X_A,\,.\,)=-{\rm d}f_A.
 \end{equation}
Thus, the vector field $X_A$ belongs to ${\mathfrak X}_{\omega}$
and   ${\bf G}\subset{\mathcal H}_{\omega}(M)$.

\smallskip

 \begin{Thm}\label{PropHomega}  Let us assume that $G$, $H$, $\Psi$  and ${\mathfrak l}$ satisfy the hypotheses
 of Theorem \ref{pi1G} and ${\rm dim}\,\Psi=1$.  If
 for each $Z\in {\mathfrak X}_{\omega}$ it is
possible to choose a function $J_Z$ on $M$, such that:
 (a) $\omega(Z,\,.\,)= -{\rm d}J_Z$;
(b) the map
 $Z\mapsto J_Z$ is continuous;
 c) $J_{X_A}=f_A$, for all $A\in{\mathfrak g}$.
  Then
  $\#\,\pi_1({\mathcal
 H}_{\omega}(M))\geq\#\,\Psi(Z(G))$.
 \end{Thm}

{\it Proof.}
 For
$Z\in{\mathfrak X}_{\omega}$, we set
 \begin{equation}\label{a(Z):}
 a(Z):[g,\alpha]\in{\mathcal
P}\longmapsto J_Z(gH)\in{\mathbb C}.
 \end{equation}
 On the other hand,
 %by (\ref{bfK}) and (\ref{omega(XA,XB}), the curvature ${\bf K}$ projects on
 %$M$ precisely the $2$-form $\omega$.
 the relation (a) is
 equivalent to $Da(Z)=-{\bf K}(Z^{\sharp},\,.\,)$. Hence,
 ${\mathfrak X}_{\omega}$
  %satisfies the conditions imposed in
 %Definition 1; i.e. it
 admits a lift to ${\mathcal P}$. The theorem follows from
Theorem \ref{pi1G}.

  \qed

\smallskip

 When  $G$ is a  complex semisimple Lie group and $P$ a parabolic
 subgroup, we put ${\mathcal F}$ for the flag manifold $G/P$. Then
${\mathcal F}\simeq U_{\mathbb R}/H$, where $U_{\mathbb R}$ is a
real compact form of $G$ and $H=P\cap U_{\mathbb R}$. Since
$U_{\mathbb R}$ and $H$ are compact groups, they are unimodular.
Hence, ${\mathcal F}= U_{\mathbb R}/H$ has a nonzero left
invariant Borel measure $d\mu$ (see \cite [Theorem 8.36]{Kn}).
  The subgroup $H$ contains a maximal torus of $U_{\mathbb R}$, so, $Z(G)\subset H$; and
  by means of the system of roots associated to this torus, it is possible to define a
  complement ${\mathfrak l}$
  of ${\mathfrak h}$ in ${\mathfrak u}_{\mathbb R}$, invariant under ${\rm Ad}(H)$.

Given $\Psi$ a representation of $H$ of dimension $1$ (as $H$ is
compact, we can assume that $\Psi$ is unitary), we have
 the respective presymplectic form $\omega$ and ${\mathfrak
 X}_{\omega}$.
  In this case,
 a choice of the functions $J_Z$ which satisfy the
conditions (a)-(c) in the statement of Theorem \ref{PropHomega}
can be carried out as follows.  Given $Z\in{\mathfrak
X}_{\omega}$, we fix $J_Z$ imposing the following normalization
condition
 \begin{equation}\label{normJZ}
 \int_{\mathcal F}J_Zd\mu=0.
 \end{equation}
 On the
other hand, denoting with $l_g$ the left translation on
$U_{\mathbb R}$ defined by $g\in U_{\mathbb R}$, for
$A\in{\mathfrak u}_{\mathbb R}$
$$I(A):=\int_{\mathcal F}f_Ad\mu=\int_{\mathcal F}f_A\circ l_g\,d\mu=\int_{\mathcal F} f_{g^{-1}\cdot
A}\,d\mu=I(g^{-1}\cdot A).$$
 Hence $I[A,B]=0$, for all $B\in{\mathfrak u}_{\mathbb R}$. Thus, $I$ vanishes on the
 commutator ideal
 $[{\mathfrak u}_{\mathbb R},\,{\mathfrak u}_{\mathbb R}]$, that coincides with ${\mathfrak
 u}_{\mathbb R}$. That is, $f_A$ satisfies the normalization condition
 (\ref{normJZ}).

 From the preceding result together with Theorem \ref{PropHomega}, it follows the
 following theorem:
 \begin{Thm}\label{Thmcartansubgroup}
 Let ${\mathcal F}$ be the flag manifold $U_{\mathbb  R}/H$, where $U_{\mathbb R}$ is a compact real form of the
 complex semisimple Lie group $G$.   Let $\Psi$ be a $1$-dimensional   representation
 of $H$ and $\omega$ the presymplectic form determined by $\Psi$ and ${\mathfrak l}$. If   ${\mathcal G}$ is any subgroup of ${\mathcal
 H}_{\omega}({\mathcal F})$ satisfying (\ref{aligned}) (ii), then
$$\#\,\pi_1({\mathcal
 G})\geq\#\,\Psi(Z(U_{\mathbb R})).$$
 \end{Thm}

 The linear map $\Psi'$ can be
extended to an element $\eta$ of ${\mathfrak u}_{\mathbb R}^*$  by
putting $\Psi'|_{\mathfrak l}=0$. If the stabilizer of $\eta$   with respect
to the coadjoint action  is precisely $H$, then
 ${\mathcal F}$ is the coadjoint orbit associated to  $\eta$ and
 $\omega$ defines a symplectic structure on
${\mathcal F}$, namely $-2\pi i\varpi$, where $\varpi$ is the
Kirillov structure. Under this hypothesis, ${\mathcal
H}_{\omega}({\mathcal F})$ is the Hamiltonian group ${\rm
Ham}({\mathcal F},\,\varpi)$ of the coadjoint orbit of $\eta$. In
this case, as invariant measure $d\mu$ we can take the one
determined by $\omega^n$, with $2n={\rm dim}\,{\mathcal F}$.

\begin{Cor}\label{CorHam}
Under the hypotheses of Theorem   \ref{Thmcartansubgroup}, if $H$
is the stabilizer of $\eta$, then
 $$\#\,\pi_1({\rm Ham}({\mathcal
 F},\,\varpi))\geq\#\,\Psi(Z(U_{\mathbb R})).$$
\end{Cor}

\smallskip

\noindent
 {\bf Remark.}  For $U_{\mathbb R}={\rm SU}(2)$,  and
$H={\rm U}(1)\subset {\rm SU}(2)$, the corresponding flag manifold
is ${\mathbb C}P^1$. For $[z_0:z_1]\in {\mathbb C}P^1$ with
$z_0\ne 0$, we put $x+iy=z_1/z_0$. It is straightforward to verify
that the vector fields $X_C$ and $X_D$ defined by the matrices of
$\mathfrak{su}(2)$
$$C=\begin{pmatrix}0& c\\ -c& 0\end{pmatrix},\;\;\;D=\begin{pmatrix} 0 & di \\ di &
0 \end{pmatrix}.$$
 take at the point $(x=0,\,y=0)$ the values
 $X_C=-c\,\frac{\partial}{\partial x}$, $X_D=d\,\frac{\partial}{\partial y}.$
 % We denote by $\omega$
% the Fubini-Study form on ${\mathbb C}P^1$, then
% $\omega_{[1:\,0]}(X_C,\,X_D)=-cd/\pi$.

% On the other hand,
  Let $\Psi$ be the character of $H$ defined by
  $\Psi({\rm diagonal}({\rm e}^{ai},\,{\rm e}^{-ai}))={\rm e}^{ai}$. Then
 $$\varpi_{[1:\,0]}(X_C,\,X_D)=\frac{i}{2\pi}\Psi'([C,D])=-cd/\pi.$$
  By Corollary \ref{CorHam},
  %By Theorem \ref{ThmpiHam},
  \begin{equation}\label{CP1}
  \#\,\pi_1({\rm
   Ham}({\mathbb C}P^1,\,\varpi))\geq 2.
   \end{equation}

We denote by $\sigma$
 the Fubini-Study form on ${\mathbb C}P^1$, then
 $\sigma_{[1:\,0]}(\partial_x,\,\partial_y)=1/\pi$.
 By the invariance
  of $\sigma$ and $\varpi$ under the action of ${\rm SU}(2)$, we
  conclude that $\varpi=\sigma$.
  % So, the Hamiltonian groups of $({\mathbb C}P^1,\,\omega)$ and $({\mathbb
  % C}P^1,\,\varpi)$ are isomorphic.
 So, (\ref{CP1})   is consistent
   with the fact that  $\pi_1({\rm
   Ham}({\mathbb C}P^1,\,\sigma))$ is isomorphic to ${\mathbb
   Z}/2{\mathbb Z}$ \cite[page 52]{lP01}.

\medskip

\noindent
 {\bf Remark.} The result stated in Corollary \ref{CorHam}
 %\ref{CorFlag}
 can be deduced using the invariancy of the action integral under loop deformations.
 More precisely, let $\psi$ be the loop in ${\rm Ham}({\mathcal
 F,\,\varpi})$
 %${\mathcal H}_{\omega}({\mathcal F})$
 at the identity,
generated by a time-dependent vector field $Z_t$ on ${\mathcal
F}$.
Here, we fix the Hamiltonian
$J_{Z_t}$ imposing the   condition
% We set $J_t$ for
%the normalized Hamiltonian $J_{Z_t}$; i.e.
 \begin{equation}\label{normaJ}
 \int_{\mathcal F} J_{Z_t}\varpi^{n}=0,
  \end{equation}
 for all $t$.
As in Subsection \ref{SubSeccAction}, we put $x_t$ for the closed
evaluation curve $\psi_t( x)$, where $x$ is a given point in
${\mathcal F}$. Since ${\mathcal F}$ is simply connected
\cite[page 33]{G-L-S}, there is a $2$-chain $C$ in ${\mathcal F}$
whose boundary is the curve $\{x_t\}$.
 Thus, the right hand side of (\ref{Theta}) is the exponential of the action
integral around the loop $\psi$ \cite{V, We}.

Let $\kappa^s$ be a deformation of $\psi$, as in Subsection
\ref{SubSeccAction}. Since  $\kappa^s_t$ is a
$\varpi$-symplectomorphism  and
$$\int_{\mathcal F} J_{Z^s_t}\varpi^{n}=0,$$
then
 \begin{equation}\label{norDotJ}
 \int_{\mathcal F}(\Dot J_t \circ\psi_t)\varpi^n=0,
  \end{equation}
 where $\Dot J_t$ is the function defined in (\ref{DotJt}).

By Proposition \ref{dThetads},   the function
 $$m:=\int_0^1(\Dot J_t\circ\psi_t)\,dt:x\in {\mathcal F}\longmapsto
 \int_0^1 \Dot J_t(\psi_t(x))\,dt,$$
 is constant. Hence, by (\ref{norDotJ})
 $$m\int_{\mathcal F}\varpi^n=\int_{\mathcal F}\Big(\int_0^1(\Dot J\circ\psi_t)\,dt
 \Big)\varpi^n=\int_0^1dt\Big(\int_{\mathcal F}(\Dot J_t\circ\psi_t)\varpi^n \Big)=0.$$
Thus,
$$m=\frac{1}{\Theta(\psi)}\frac{{\rm d}\Theta^s}{{\rm
d}s}\bigg|_{s=0}=0,$$
% Now, Corollary \ref{CorFlag} can be deduced using the fact that the
and we obtain the known fact that the  action integral is
invariant under deformations of the loop; that is, $\Theta$
factors through the homotopy group, giving rise to the map
(\ref{Theta:pi_1}). From (\ref{F1(p)=}), it follows that
$\Theta(\psi\star\tilde\psi)=\Theta(\psi)\Theta(\tilde\psi)$,
where $\star$ denotes the path
  product. Thus, (\ref{Theta:pi_1}) is  a group homomorphism.

Now, using this invariance of $\Theta$, we can deduce Corollary
\ref{CorHam}. Given $g_1\in Z(U_{\mathbb R})$, let
 $A\in{\mathfrak h}$, such that ${\rm e}^A=g_1$. We set $g_t={\rm
 e}^{tA}$; this path defines the corresponding loop $\varphi$ in
 the Hamiltonian group and the respective Hamiltonian function is
 $J_t=f_{tA}=t\Psi'(A)$.
  The curve $\varphi_t(\bar e)$ reduces to
 $\{\bar e\}$ and the corresponding value of $\Theta$ given by (\ref{Theta}) is
 $\Theta(\varphi)=\Psi(g_1)$. Thus, the loops $\varphi$ associated to two
 elements in $Z(U_{\mathbb R})$ on which $\Psi$ takes different
 values, determine distinct elements in $\pi_1({\rm Ham}({\mathcal
 F},\,\varpi)).$

\medskip

%%%%%%%%%%%%%%%%%%%%%%%%%%%%%%%%%%%%%%%%%%%%%%%%%%%%%%%%%%%%%%%%%%%%%%%%%%%%%%%%%%%%%%%%%%%%%%%%%%%%%%%%%%%
%%%%%%%%%%%%%%%%%%%%%%%%%%%%%%%%%%%%%%%%%%%%%%%%%%%%%%%%%%%%%%%%%%%%%%%%%%%%%%%%%%%%%%%%%%%%%%%%%%%%%%%%%%%%%%

\section{Toric manifolds}\label{SectToric}

In this section,  we will consider symplectic toric manifolds and the group ${\mathcal G}$ will be the respective Hamiltonian group.

Let $T$ be the torus $({\rm
 U}(1))^n$, and  $\Delta$ a Delzant polytope in ${\mathfrak t}^*$
 \cite{Del}, whose mass center is denoted by ${\rm Cm}(\Delta)$.
 We denote by   $(M_{\Delta},\,\omega_{\Delta})$  the toric
$2n$-manifold, defined by the
  polytope $\Delta$ \cite{Gui}. We put  $\Phi$ for the corresponding moment map
 $\Phi:M\to {\mathfrak t}^*$.

%We will consider the group ${\mathcal G}={\rm Ham}(M_{\Delta},\,
%\omega_{\Delta})$.
 Since $M_{\Delta}$ is simply connected, ${\rm
Lie}({\rm Ham}(M_{\Delta},\, \omega_{\Delta}))$  consists of the
vector fields $Z$ such that ${\rm d}\iota_Z\omega_{\Delta}=0$. For such a $Z$,
we denote by $f_Z$ the corresponding mean normalized Hamiltonian;
i.e.
 \begin{equation}\label{dfZ}
 {\rm d}f_Z=-\iota_Z\omega_{\Delta},\;\;\;
 \int_{M_{\Delta}}f_Z(\omega_{\Delta})^{n}=0.
 \end{equation}

 Let ${\bf b}$ be a vector  of ${\mathfrak
  t}$.
  We set $X_{\bf b}$ for the vector field on
$M_{\Delta}$ generated by ${\bf b}$ through the  $T$-action;
$X_{\bf b}$ is a Hamiltonian vector field and
$${\rm d}\langle\Phi,\,{\bf b}\rangle=-\omega_{\Delta}(X_{\bf
b},\,.\,).$$
 The respective
 normalized Hamiltonian
 %for the circle action generated by ${\bf b}$
  is the function
  \begin{equation}\label{fbfb}
  f_{\bf b}:= \langle \Phi,\,{\bf b}\rangle-\langle{\rm Cm}(\Delta),\,{\bf
  b}\rangle,
   \end{equation}
    where
 \begin{equation}\label{Cm}
\langle \text{Cm}(\Delta),\,{\bf b}\rangle= \frac{\int_M
\langle\Phi,\,{\bf b}\rangle\,(\omega_{\Delta})^n}{\int_M
(\omega_{\Delta})^n}.
\end{equation}

\smallskip

%%%%%%%%%%%%%%%%%%%%%%%%%%%%%%%%%%%%%%%%%%%%%%%%%%%%%%%%%%%%%%%%%%%%%%%%%%%%%%%%%%%%%%%%%%%%%%%%%%%%%%%%%%%%%%%%%%%%%%%%

\subsection{Quantizable manifolds}\label{SubsecQuantizable}
 Identifying ${\mathfrak t}^*$ with ${\mathbb R}^n$, we
 can assume that $*:=(0,\dots,0)$ is a vertex of $\Delta$. Under this assumption,
  the symplectic manifold $(M_{\Delta},\,\omega_{\Delta})$ is quantizable
 iff the vertices of $\Delta$ have integer coordinates.

 In this subsection,
 we assume that $(M_{\Delta},\,\omega_{\Delta})$ {\em is quantizable}. Then
there exists a $U(1)$
 principal bundle  ${\mathcal L} \stackrel{{\rm pr}}{\longrightarrow}  M_{\Delta}$ with a connection such
 that the curvature ${\bf K}$ projects a $2$-form $\omega$ on
 $M_{\Delta}$, with
 \begin{equation}\label{i2pi}
  \frac{i}{2\pi}\omega=\omega_{\Delta}.
  \end{equation}

We will show that ${\mathfrak X}$, the Lie algebra of the ${\rm
Ham}(M_{\Delta},\,\omega_{\Delta})$, ``admits a lift" to
${\mathcal L}$, in the sense that it satisfies properties similar
to the ones stated in Definition 1.
 For
 $Z\in{\mathfrak X}$,
 %{\rm Lie}(\mathcal G)$,
   we denote by $a(Z)$ the map
 $$a(Z):p\in{\mathcal L}\longmapsto -2\pi i f_Z({\rm pr}(p))\in i{\mathbb
 R}.$$
 For this map the following conditions hold:

 \noindent
 (i) $a(Z)(p\beta)=a(Z)$, for any $\beta\in U(1)$.

 \noindent
 (ii) From (\ref{i2pi}) together with (\ref{dfZ}), it follows
  $$Da(Z)=-{\bf K}(Z^{\sharp},\,.\,),$$
  where $Z^{\sharp}$ is the horizontal lift of $Z$.

\noindent
  (iii) The continuity of $Z \mapsto a(Z)$ is guarantied by   conditions (\ref{dfZ}), which determine  $f_Z$ uniquely.

\noindent
 (iv) Obviously,
  $a(X_{\bf b})(p)=-2\pi i f_{\bf b}({\rm pr}(p)).$

  % Therefore the
  %algebra ${\mathfrak X}$ of the Hamiltonian vector fields   $Z$
  %satisfies properties   analogous to de ones stated in Definition
 %1.

  \smallskip

 As in Section \ref{SecX}, we define
  $U(Z)=Z^{\sharp}+V_{a(Z)}$, where $Z^{\sharp}$ is the horizontal lift of $Z$ to ${\mathcal L}$.
   By the way,  the vector field $U(Z)$ is the operator assigned to the
   function $f_Z$ in Geometric Quantization
  (see \cite[formula (3.29)]{Sn}).
Propositions \ref{Proprighttrans}
  and
  \ref{ProprighttransOmega} in Section \ref{SecX} have an immediate translation to the
  present situation.

   If $\psi$ is the loop in ${\rm
  Ham}(M_{\Delta},\,\omega_{\Delta})$ at the identity, generated by the family $Z_t\in{\mathfrak
 X}$,   then we
  have the corresponding isotopy $F_t$ in ${\mathcal L}$ defined by the time-dependent vector field $U(Z_t)$. As in
  Proposition \ref{fconstant2}, $F_1$ is the multiplication by a constant, say
   $F_1(p)=p\Lambda(\psi)$.

   Taking into account that $H_1(M_{\Delta},\,{\mathbb Z})=0$, in the same way as in Proposition
 \ref{PropTheta}, we have:
 \begin{Prop}\label{kappa}
Let $\psi$ be the loop in ${\rm
Ham}(M_{\Delta},\,\omega_{\Delta})$ generated by the family
$Z_t\in{\mathfrak
 X}$    and $x$ an arbitrary point of
$M_{\Delta}$, then
 $$\Lambda(\psi)={\rm exp}\,\Big(2\pi
 i\int_C\omega_{\Delta}-2\pi
 i\int_0^1f_{Z_t}(x_t)\,dt\Big),$$
 where $C$ is any $2$-chain whose boundary is the curve $\{x_t=\psi_t(x)\}$.
 \end{Prop}

As in the last Remark of Subsection \ref{SubSeccHVF}, $\Lambda$
defines a group homomorphism
\begin{equation}\label{Lambda:pi_1}
\Lambda:\pi_1({\rm
Ham}(M_{\Delta},\,\omega_{\Delta}))\longrightarrow U(1).
\end{equation}

 Let ${\bf b}$ be a   vector in the integer lattice of ${\mathfrak t}$;
 it defines a loop $\varphi_{\bf b}$ in ${\rm
 Ham}(M_{\Delta},\,\omega_{\Delta})$. The vertices of
 $\Delta$ are the fixed points of the $T$-action, so $\varphi_{\bf
 b}(t)(*)=*$ for all $t$. From (\ref{fbfb}), $f_{\bf b}(*)=-\langle {\rm Cm}(\Delta),\,{\bf b}\rangle$, and
  Proposition \ref{kappa} applied to the point $*$ gives
   \begin{equation}\label{Lambdavafrphi_b}
   \Lambda(\varphi_{\bf b})={\rm exp}\,\big(2\pi i  \langle{\rm Cm}(\Delta),\,{\bf b}\rangle\big).
   \end{equation}
By (\ref{Lambda:pi_1}),  we can state the following theorem:
\begin{Thm}\label{ThmCM(Delta)}
Let ${\bf b}$, ${\bf\tilde b}$ be  vectors in the integer lattice
of ${\mathfrak t}$. If
 \begin{equation}\label{formu}
\langle{\rm Cm}(\Delta),\,{\bf b}-{\bf\tilde
b}\rangle\notin{\mathbb Z},
 \end{equation}
then
$$[\varphi_{\bf b}]\ne[\varphi_{\bf\tilde b}]\in\pi_1({\rm Ham}(M_{\Delta},\,\omega_{\Delta})).$$
 \end{Thm}

As $(M_{\Delta},\,\omega_{\Delta})$ is quantizable, the coordinates of ${\rm Cm}(\Delta)$ are rational numbers, say
 ${\rm Cm}(\Delta)=(\frac{r_1}{r},\dots,\frac{r_n}{r})$. From Theorem \ref{ThmCM(Delta)}, we deduce the
 corollary:
 \begin{Cor}\label{Corgcd}
 If ${\rm g.c.d.}(r_1,\dots,r_n)=1$, then there exists ${\varphi_{{\bf b}_j}}$, with $j=1,\dots,|r|$, such that for $i\ne j$
 $$[{\varphi_{{\bf b}_i}}]\ne[{\varphi_{{\bf b}_j}}]\in \pi_1({\rm Ham}(M_{\Delta},\,\omega_{\Delta})).$$
 \end{Cor}

\smallskip

\noindent {\bf Remark.} The values taken by the map $\Lambda$
admit an interpretation in terms of the integration of
$T$-equivariant forms on $M_{\Delta}$. As $T$ is abelian, the map
$f=\Phi-{\rm Cm}(\Delta):M\to {\mathfrak t}^*$ is $T$-equivariant.
By (\ref{dfZ})(i), $\omega_{\Delta}-f$ is an $T$-equivariantly
closed $2$-form on $M_{\Delta}$ \cite[page 111]{G-S}.
 Given $q\in M_{\Delta}$ a fixed point for the $T$-action; let $m_{q,1},\dots, m_{q,n}$ denote the weights of the isotropy
 representation of $T$ on the tangent space $T_qM_{\Delta}.$
 The
localization formula of equivariant cohomology \cite{B-G-V} applied to the form $\alpha={\rm
exp}\big(2\pi i(\omega-f) \big)$ gives, for a lattice vector ${\bf
b}$,
$$\int_M \alpha({\bf b})=\sum_q\frac{{\rm exp}(-2\pi if_{\bf
b}(q))}{\prod_k m_{q,k}({\bf b})},$$
 where $q$ runs over the set of fixed points of the $T$-action (assumed that
 all $ m_{q,k}({\bf b})$ are non zero).

 Since the manifold is quantizable and ${\bf b}$ is in the integer
 lattice, for two fixed points $q$ and $q'$ of the $T$-action,
 $f_{\bf b}(q)-f_{\bf b}(q')\in{\mathbb Z}$. So, by
 (\ref{Lambdavafrphi_b}),
 $$\int_M \alpha({\bf b})=\Lambda(\varphi_{\bf b})\sum_q\Big(\prod_k m_{q,k}({\bf b})
 \Big)^{-1}.$$

By means of the localization theorem for $S^1$-actions, it is possible to deduce sufficient conditions for
two $S^1$-actions   $\varphi_{\bf b}$ and  $\varphi_{\bf\tilde b}$ not to be homotopically equivalent through
a homotopy consisting of {\it  $S^1$-actions} (see \cite[Theorem 5]{V08}).
% However,
The conclusion of Theorem \ref{ThmCM(Delta)}
is obviously stronger, although it is only applicable to quantizable manifolds.

\medskip
\noindent
 {\bf Examples.}  % In the examples of Section \ref{SectToric}, we consider some
    % There are examples of
    For some toric manifolds $M_{\Delta}$,
     we  proved in \cite{V2b}   the existence of infinite cyclic subgroups in $\pi_1({\rm Ham}(M_{\Delta}))$, generated by lattice
     vectors ${\bf b}$ such that $(\Delta,\,{\bf b})$ {\em  is not} a mass linear pair. Below, in the Examples 2 and 3,
      we consider some  of those
     % examples
      manifolds and show that there are  lattice
     vectors ${\bf b}$ that define  homotopically non trivial loops in ${\rm Ham}(M_{\Delta})$, although $(\Delta,\,{\bf b})$
     are  mass linear pairs \cite{M-T3}.

\smallskip

 \noindent
 {\it Example 1: The projective space ${\mathbb C}P^n$.} The manifold $M_{\Delta}$ associated to the standard simplex
$$\Delta=\Big\{(x_1,\dots,x_n)\in{\mathbb R}^n\,|\,\sum_ix_i\leq
1,\,0\leq x_i\Big\},$$ is ${\mathbb C}P^n$. The mass center ${\rm
Cm}(\Delta)=(\frac{1}{n+1},\dots,\frac{1}{n+1})$. Denoting ${\bf
b_j}=(j,0,\dots,0)$, then $\langle{\rm Cm}(\Delta),\,{\bf
b}_j\rangle=\frac{j}{n+1}$. By Theorem \ref{ThmCM(Delta)}, the
loops $\varphi_{{\bf b}_i}$, for $i=0,1\dots,n$ determine $n+1$
different elements in $\pi_1({\rm
Ham}(M_{\Delta},\,\omega_{\Delta}))$. This result coincides with
the one stated in \cite[Theorem 5]{V2}, which was deduced by
calculations of Maslov indices. Furthermore, that lower bound for
the cardinal of $\pi_1({\rm Ham}(M_{\Delta},\,\omega_{\Delta}))$
is consistent with the following two facts: $\pi_1({\rm
Ham}({\mathbb C}P^1))={\mathbb Z}/2{\mathbb Z}$, and ${\rm
Ham}({\mathbb C}P^2)$ has the homotopy type of ${\rm PU}(3)$
\cite{Gr}.

\smallskip
\noindent
 {\it Example 2: Hirzebruch surfaces.} Given $r\in{\mathbb Z}_{>0}$ and
 $\tau,\lambda\in{\mathbb R}_{>0}$ with $\sigma:=\tau-r\lambda>0$.
 The polytope $\Delta$ in ${\mathbb R}^2$ defined by the
 vertices
 $$(0,\,0),\; (0,\,\lambda),\; (\tau,\,0),\; (\sigma,\,\lambda)$$
  determines  a Hirzebruch surface  $M_{\Delta}$ equipped with the corresponding symplectic structure.
The mass center of $\Delta$ is
 \begin{equation}\label{Cmblowup}{\rm
Cm}(\Delta)=\Big(\frac{3\tau^2-3r\tau\lambda+r^2\lambda^2}{3(2\tau-r\lambda)},\;
\frac{3\lambda\tau-2r\lambda^2}{3(2\tau-r\lambda)} \Big).
\end{equation}

In \cite[Corollary 4.2]{V2b} and using the homomorphism
(\ref{Idefinition}), we proved that $\varphi_{\bf b}$ generates an
infinite cyclic subgroup in $\pi_1({\rm
Ham}(M_{\Delta},\,\omega_{\Delta}))$, if ${\bf b}=(b_1,\,b_2)$
satisfies $rb_1\ne 2b_2$. The point was that the homomorphism
(\ref{Idefinition})
%$I(\varphi_{\bf b}
vanishes on $[{\varphi}_{\bf b}]$ precisely when $rb_1= 2b_2$; in
other words, when $(\Delta,\,{\bf b})$ is a mass linear pair.

 However, when $(M_{\Delta},\,\omega_{\Delta})$ is quantizable, if $rb_1=
2b_2$, then  the map (\ref{Lambda:pi_1}) on $\varphi_{\bf b}$  takes
the value
$$\Lambda(\varphi_{\bf b})={\rm exp}\big(2\pi i\langle {\rm Cm}(\Delta),\,{\bf b}
\rangle\big)={\rm exp}(\pi ib_1\tau).$$
 Thus, if $b_1\tau$ is odd and $rb_1$ even, then the lattice vector ${\bf
 b}=(b_1,\,\frac{rb_1}{2})$ determines a loop   homotopically
non trivial in the Hamiltonian group.

%To summarize,
 In the context of the toric quantizable manifolds, the homomorphisms $I$ and $\Lambda$ have
``complementary properties", in the following sense:
% when one considers its actions on $\pi_1({\rm
%Ham}(M_{\Delta},\,\omega_{\Delta}))$.
% on quantizable manifolds.
  The map (\ref{Idefinition})
 can not detect nontrivial elements of finite order in $\pi_1({\rm Ham}(M_{\Delta},\,\omega_{\Delta}))$, because
  it is an ${\mathbb R}$-valued group
homomorphism. By contrast, $\Lambda$ can not distinguish elements
$[\varphi_{\bf b}]$ of finite order from those of order infinite. In fact, as ${\rm Cm}(\Delta)$ has rational coordinates, for any lattice vector ${\bf b}$
 there is an integer $m$, such that $\Lambda([\varphi_{\bf
b}]^m)=1$ (see (\ref{Lambdavafrphi_b})).

\smallskip
\noindent
 {\it Example 3: One point blow-up of ${\mathbb C}P^n$.} Given $\tau,\lambda\in{\mathbb Z}_{>0}$ with $\sigma:=\tau-\lambda>0$, we denote by
 $\Delta$ the following truncated simplex
 $$\Delta=\Big\{(x_1,\dots,x_n)\in{\mathbb R}^n\,|\, \sum_i x_i\leq\tau,\;0\leq x_i,\; x_n\leq \lambda\Big\}.$$
  The corresponding toric manifold $M_{\Delta}$ is the one point blow-up of  ${\mathbb C}P^n$.

  The mass center ${\rm Cm}(\Delta)$ is the point
  $$ {\rm Cm}(\Delta)=\frac{1}{\tau^n-\sigma^n}\Big(\Big(\frac{\tau^{n+1}-\sigma^{n+1}}{n+1}\Big)w-\lambda\sigma^n e_n\Big),$$
 where $w=(1,\dots, 1)$ and $e_n=(0,\dots,0,1)$. In \cite{V2b}, we proved that the loop $\varphi_{\bf b}$,
 defined by  a lattice
  vector ${\bf b}=(b_1,\dots, b_n)$, with $n b_n\ne\sum_{i=1}^{n-1}b_i$,  determines  an element  $[\varphi_{\bf b}]$ in
  $\pi_1({\rm Ham}(M_{\Delta},\,\omega_{\Delta}))$ of infinite
  order. This property  is consequence of the non-vanishing of the
  group homomorphism $I$ on $[{\varphi}_{\bf b}]$.

  When $(\Delta,\,{\bf b})$ is a linear pair; i. e., $n b_n=\sum_{i=1}^{n-1}b_i$,
  \begin{equation}\label{langerm}
  \langle{\rm Cm}(\Delta),\,{\bf b}\rangle=
  \Big(\frac{\tau^{n+1}-\sigma^{n+1}-\lambda\sigma^n}{\tau^n-\sigma^n}\Big)b_n.
  \end{equation}
Hence, if (\ref{langerm}) is no integer, then ${\varphi_{\bf b}}$
is not contractible.

%%%%%%%%%%%%%%%%%%%%%%%%%%%%%%%%%%%%%%%%%%%%%%%%%%%%%%%%%%%%%%%%%%%%%%%%%%%%%%%%%%%%%%%%%%%%%%%%%%%%%%%%%%%%%%%%%%%%%%%%%

\smallskip

\subsection{Some non quantizable manifolds}\label{SubsecNonQuantizable}
Let us assume that $\Delta$ has $d$ facets. Denoting with ${\bf n}_j\in{\mathfrak t}$ the normalized conormals to the facets, then the
correspondence $e_j\mapsto {\bf n}_j$, where $\{e_j\}$ is the standard basis of ${\mathbb R}^d$ induces a group homomorphism
$${\mathbb R}^d/{\mathbb Z}^d\longrightarrow T={\mathbb R}^n/{\mathbb Z}^n,$$
whose kernel will be denoted by $N$ (see \cite[Chapter 1]{Gui}).
The standard action of ${\mathbb R}^d/{\mathbb Z}^d$ on ${\mathbb
C}^d$ gives rise to an action of $N$ on that space, which
preserves the splitting
$$\bigoplus_{j=1}^d {\mathbb C}e_j={\mathbb C}^d.$$
  Let ${\mathcal L}_j$
denote the line bundle over $M_{\Delta}$ defined by the
representation of $N$ on ${\mathbb C}e_j$.
%\subset {\mathbb C}^d$.
We denote by $c_j$ the Chern class of ${\mathcal L}_j$. The
classes $c_1,\dots,c_d$ generate the ring
$H^*(M_{\Delta},\,{\mathbb C})$ (see \cite[Proposition
9.8.7]{G-S}).

Let us assume that the cohomology class $[\omega_{\Delta}]$ of the symplectic structure satisfies
\begin{equation}\label{[omegaDelta]}
[\omega_{\Delta}]=r\sum_{j=1}^d n_jc_j,
 \end{equation}
 with $r\in{\mathbb R}$ and $n_j\in{\mathbb Z}.$ Obviously, the sum in (\ref{[omegaDelta]}) is the Chern class of the line bundle
 $${\mathcal L}:=\bigotimes_{j=1}^d{\mathcal L}_j^{\otimes n_j}.$$

 The manifold $(M_{\Delta},\,\omega_{\Delta}/r)$ is a quantizable and there exists a connection on ${\mathcal L}$
  such that its curvature ${\bf K}$ projects on $M_{\Delta}$ the form $-2\pi i\omega_{\Delta}/r$.

  Let ${\mathfrak X}$ denote the algebra of Hamiltonian vector fields on $(M_{\Delta},\,\omega_{\Delta})$. As at the beginning  of the section, $f_Z$
  will be normalized Hamiltonian function associated to $Z\in{\mathfrak X}$; i.e. $f_Z$ satisfies (\ref{dfZ}).  The map
  $$a(Z): p\in{\mathcal L}\longmapsto -\frac{2\pi i}{r}f_Z({\rm pr}(p))\in{\mathbb C},$$
   is  obviously  invariant under the right translations in ${\mathcal L}$ and $D a(Z)=-{\bf K}(Z^{\sharp},\,.\,)$.

  Given $\psi$ a loop in ${\rm Ham}(M_{\Delta},\,\omega_{\Delta})$ at the identity generated by a family $Z_t$
  of vector fields, as in Subsection \ref{SubsecQuantizable}, the corresponding gauge transformation $F_1$
  is the multiplication by a constant $\Lambda(\psi)$, with
   $$\Lambda(\psi)={\rm exp}\,\Big(\frac{2\pi
 i}{r}\int_C\omega_{\Delta}-\frac{2\pi
 i}{r}\int_0^1f_{Z_t}(x_t)\,dt\Big),$$
 where $C$ is any $2$-chain whose boundary is $\{x_t=\psi_t(x)\}$.
 Hence, we have the following proposition:
 \begin{Prop}\label{Propnonquantiz}
 Assumed that  $[\omega_{\Delta}]$ satisfies (\ref{[omegaDelta]}). If  ${\bf b}$ and  ${\bf\tilde b}$ are lattice vectors
 such that
 \begin{equation}\label{forPropnonq}
 \langle {\rm Cm}(\Delta),\,{\bf b}-{\bf \tilde b}\rangle\notin r{\mathbb Z},
  \end{equation}
 then
 $$[\varphi_{\bf b}]\ne[\varphi_{\bf\tilde b}]\in\pi_1({\rm
 Ham}(M_{\Delta},\,\omega_{\Delta})).$$
  \end{Prop}

\smallskip
\noindent
 {\bf Remark.} Let $\Delta'$ denote the polytope
 $\frac{1}{r}\Delta$. The manifold
 $(M_{\Delta'},\,\omega_{\Delta'})$ is quantizable and the
 hypothesis (\ref{formu})
  % condition $\langle{\rm Cm}(\Delta'),\,{\bf b}-{\bf \tilde
% b}\rangle\notin{\mathbb Z}$
 of Theorem \ref{ThmCM(Delta)} applied to  $(M_{\Delta'},\,\omega_{\Delta'})$
  is equivalent to (\ref{forPropnonq}).

\bigskip

%%%%%%%%%%%%%%%%%%%%%%%%%%%%%%%%%%%%%%%%%%%%%%%%%%%%%%%%%%%%%%%%%%%%%%%%%%%%%%%%%%%%%%%%%%%%%%%%%%%%%%%%%%%%%%%%%%%%%%%%%%%%%%%%%%%
%%%%%%%%%%%%%%%%%%%%%%%%%%%%%%%%%%%%%%%%%%%%%%%%%%%%%%%%%%%%%%%%%%%%%%%%%%%%%%%%%%%%%%%%%%%%%%%%%%%%%%%%%%%%%%%%%%%%%%%%%%%%%%%%%%%

\end{document}